\documentclass[10pt]{amsart}
\pdfoutput=1
\usepackage{amsmath, amsthm, amssymb,slashed}
\usepackage{ifpdf}
\usepackage[pdftex]{graphicx}
\usepackage{tikz}
\usetikzlibrary{matrix,arrows,calc}
\usepackage[pdftex,plainpages=false,hypertexnames=false,pdfpagelabels]{hyperref}
 \theoremstyle{plain}
 \newtheorem{thm}{Theorem}[section]
 
 \newtheorem{lem}[thm]{Lemma}
 \newtheorem{prop}[thm]{Proposition}
 \theoremstyle{definition}
 \newtheorem{defn}[thm]{Definition}
 \newtheorem{notation}[thm]{Notation}

 \newtheorem{ex}[thm]{Example}
 \theoremstyle{remark}
 \newtheorem{rmk}[thm]{Remark}

\def\beq{\begin{eqnarray}}
\def\eeq{\end{eqnarray}}

\DeclareSymbolFont{bbold}{U}{bbold}{m}{n}
\DeclareSymbolFontAlphabet{\mathbbold}{bbold}
\def\one{\mathbbold{1}}
 \newcommand{\bp}{\begin{proof}[Proof]}
 \newcommand{\ep}{\end{proof}}
 \DeclareMathOperator{\Hess}{Hess}
\DeclareMathOperator{\SM}{\underline{\sf SMfld}}

\newcommand{\sq}{\mathord{/\!\!/}}

\def\odd{{\rm odd}}
\def\ev{{\rm ev}}
\def\cl{{\rm cl}}

\def\pt{{\rm pt}}
\def\vol{{\rm vol}}
\def\partition{{Z}}
\def\R{{\mathbb{R}}}

\def\id{{{\rm id}}}
\def\proj{{\rm proj}}
\def\MF{{\rm MF}}
\def\MaF{{\rm MaF}}

\def\C{{\mathbb{C}}}
\def\Q{{\mathbb{Q}}}
\def\T{{\mathbb{T}}}
\def\Z{{\mathbb{Z}}}
\def\E{{\mathbb{E}}} 
\def\M{{\mathbb{M}}} 

\def\pt{{\rm pt}}
\def\End{\mathop{\sf End}}

\def\EB{\hbox{-{\sf EBord}}}

\def\SM{ {\underline{\sf SMfld}}}

\def\EFT{ \hbox{-{\sf EFT}}}

\def\Spin{{\rm Spin}}
\def\SO{{\rm SO}}
\def\SL{{\rm SL}}
\def\TMF{{\rm TMF}}

\def\Iso{ {\underline{\rm Iso}}}

\begin{document}

\title{Perturbative sigma models, elliptic cohomology and the Witten genus}

\author{Daniel Berwick-Evans}

\address{Department of Mathematics, University of Illinois at Urbana-Champaign}

\email{danbe@illinois.edu}

\date{\today}

\maketitle 

\begin{abstract}
We provide a differential cocycle model for elliptic cohomology with complex coefficients and use analytic methods to construct a cocycle representative for the Witten class in this language. Our motivation stems from the conjectural connection between 2-dimensional field theories and elliptic cohomology originally due to G. Segal and E. Witten. The specifics of our constructions are informed by the work of S. Stolz and P. Teichner on super Euclidean field theories and K. Costello's construction of the Witten genus using perturbative quantization. As a warm-up, we prove analogous results for supersymmetric quantum mechanics and K-theory with complex coefficients.
\end{abstract}

\setcounter{tocdepth}{1}
\tableofcontents

\section{Introduction}

In this paper we provide a differential cocycle model for elliptic cohomology with complex coefficients and use analytic methods to construct a cocycle representative for the Witten class in this language. Our motivation stems from the conjectural connection between sigma models and elliptic cohomology originally due to E. Witten~\cite{witten_ell} and G. Segal~\cite{segal_elliptic}. The specifics of our constructions are informed by the work of S. Stolz and P. Teichner on $2|1$-dimensional super Euclidean field theories~\cite{ST11} and K. Costello's construction of the Witten genus using methods of perturbative (specifically, \emph{1-loop}) quantization~\cite{costello_WG1,costello_WG2}. Indeed, our investigation began as an attempt to develop a bridge between these two languages that relate particular 2-dimensional quantum field theories to elliptic cohomology. As a warm-up case, we will prove analogous results for $1|1$-dimensional super Euclidean field theories, (twisted, Real) ${\rm K}$-theory with complex coefficients, and the $\hat{A}$-class, showing that 1-loop quantization gives rise to the familiar local index from the Atiyah--Singer theorem. 

There are two main geometric constructions in the paper we wish to highlight. The first associates to any smooth manifold~$X$ a super stack, denoted $\Phi_0^{2|1}(X)$, consisting of maps from $2|1$-dimensional tori to~$X$ that are classical solutions with zero energy for the supersymmetric sigma model with target~$X$. These are the \emph{classical vacua.} There is a natural line bundle over~$\Phi_0^{2|1}(X)$, and sections of tensor powers that are invariant under infinitesimal anti-holomorphic dilations of tori are differential cocycles for ${\rm TMF}(X)\otimes \C$, where ${\rm TMF}$ denotes the cohomology theory of topological modular forms constructed by P.~Goerss, M.~Hopkins and H.~Miller. The second main construction takes as input a Riemannian structure on $X$ and outputs a family of operators parametrized by~$\Phi_0^{2|1}(X)$, denoted~$\Delta_X^{2|1}$. 
Physically these operators encode the Hessian of the classical action of the sigma model, called \emph{kinetic operators}; more concretely, $\Delta^{2|1}_X$ is a deformation of the Laplacian on tori that incorporates the curvature 2-form of the Levi-Civita connection. The fiberwise $\zeta$-determinant of the family~$\Delta_X^{2|1}$ defines a function on~$\Phi_0^{2|1}(X)$, but generally it is not invariant under anti-holomorphic dilations of tori (so does not determine a cocycle). A rational string structure on~$X$ specifies a modification of the determinant that yields a cocycle we identify with the Witten class of~$X$ whose integral over~$X$ is the Witten genus. 

In the next subsection, we provide the necessary background to give precise statements of our main results. We spend the remainder of the section discussing physical motivation~(\S\ref{sec:introphys}) and explaining connections to previous work on the Witten genus~(\S\ref{sec:introWit}-\S\ref{intro:cost}).

\subsection{Statement of results}\label{sec:results}
Let $T^{2|1}_\Lambda=\R^{2|1}/\Lambda$ be the super torus\footnote{As usual in super geometry, one should consider \emph{families} of tori. For simplicity we will suppress this family parameter in the current discussion.} associated to a lattice~$\Lambda\subset\R^{2|1}$. Let $X$ be a smooth manifold, and $\Phi^{2|1}(X)$ denote the \emph{stack of fields} consisting of pairs $(\Lambda,\phi)$ for $\phi\colon T_\Lambda^{2|1}\to X$ a smooth map. Our primary object of study is a substack, denoted~$\Phi^{2|1}_0(X)\subset \Phi^{2|1}(X)$, of \emph{classical vacua}. This substack is defined as having maps $\phi$ factorizing through the map $T^{2|1}_\Lambda\to \R^{0|1}$ induced by the projection~$\R^{2|1}\to \R^{0|1}$. 

The stack $\Phi_0^{2|1}(\pt)$ is essentially the component of the stack of spin elliptic curves with nonbounding spin structure. As such, there is a line bundle~$\omega^{1/2}$ over $\Phi_0^{2|1}(\pt)$ whose tensor square is the Hodge bundle~$\omega$, i.e., the fiber of $\omega^{1/2}$ at $T^{2|1}_\Lambda$ is holomorphic sections of the square root of the canonical line of the underlying reduced torus~$\R^2/\Lambda\cong \C/\Lambda$, where the spin structure specifies the square root. Even tensor powers of this line bundle~$(\omega^{1/2})^{\otimes 2k}$ have as holomorphic sections weight~$k$ weak modular forms. Naturality of classical vacua permits us to pull back~$\omega^{1/2}$ to a line bundle also denoted~$\omega^{1/2}$ over~$\Phi_0^{2|1}(X)$. There are natural vector fields on $\Phi_0^{2|1}(X)$ encoding the dilation of each meridian of a torus $\R^2/\Lambda\cong S^1\times S^1$, and these decompose into holomorphic and anti-holomorphic components. We denote those sections invariant under the anti-holomorphic components by~$\mathcal{O}(\Phi_0^{2|1}(X);\omega^{\bullet/2})\subset \Gamma(\Phi_0^{2|1}(X);\omega^{\bullet/2})$ and refer to them as \emph{holomorphic sections}; see~\S\ref{sec:proofthm1}. The assignment $X\mapsto \mathcal{O}(\Phi_0^{2|1}(X);\omega^{\bullet/2})$ defines a sheaf on the site of smooth manifolds.


\begin{thm}\label{thethm}
There is a natural isomorphism of sheaves of graded algebras over~$\C$
$$
\mathcal{O}(\Phi_0^{2|1}(-);\omega^{\bullet/2})\stackrel{\sim}{\to}\bigoplus_{i+j=\bullet} \Omega^i_{\rm cl}(-)\otimes \MF^j
$$
whose target is the sheaf of closed differential forms valued in the graded ring $\MF^\bullet$ of weak modular forms. This realizes the source sheaf as a differential cocycle model for $\TMF\otimes \C$ in the sense of Hopkins--Singer~\cite{hopkins_singer}. 
\end{thm}


The part of the statement regarding Hopkins--Singer differential cocycles follows easily from the isomorphism of sheaves, as we now explain. Since ${\rm TMF}\otimes \C$ is a cohomology theory over~$\C$, it is just an ordinary cohomology theory with values in $\TMF^\bullet(\pt)\otimes \C\cong \MF^\bullet$.
By the de~Rham theorem, there is a surjective map
$$
\bigoplus_{i+j=\bullet} \Omega^i_{\rm cl}(X)\otimes \MF^j\twoheadrightarrow {\rm TMF}^\bullet(X)\otimes \C
$$
that takes a closed form to its cohomology class. 
Though the language is perhaps a bit overloaded in this easy case, this application of the de~Rham theorem along with the isomorphism of sheaves in Theorem~\ref{thethm} presents~$\mathcal{O} (\Phi_0^{2|1}(X);\omega^{\bullet/2})$ as a differential cocycle model for ${\rm TMF}(X)\otimes \C$. One might (rightfully) argue that de~Rham forms valued in modular forms are surely a simpler model. Our reason for examining this more complicated version comes from the insights it provides into the differential-geometric content of elliptic cohomology, beyond what is visible in the de~Rham model. We list a few of these features. 
\begin{enumerate}
\item The sections $\mathcal{O} (\Phi_0^{2|1}(X);\omega^{\bullet/2})$ receive a map from $2|1$-Euclidean field theories, which Stolz and Teichner~\cite{ST04,ST11} have identified as a candidate cocyle model for $\TMF(X)$ \emph{integrally}. In~\S\ref{sec:Stolz--Teichner} we conjecture that this map gives a cocycle-level description of the generalized Chern character for~TMF. 
\item Replacing~$X$ with a quotient stack $X\sq G$ provides a natural (and non-Borel) equivariant refinement of our cocycle model; see~\cite{DBE_EquivTMF}. When $G$ is a finite group, this constructs equivariant elliptic cohomology as defined by Devoto~\cite{Devoto}. 
The coefficient ring is roughly a combination of modular forms and differential forms on fixed point sets for pairs of commuting elements in~$G$, as one expects from a height~2 Hopkins--Kuhn--Ravenel character theory~\cite{HKR}. Geometrically, the pairs of commuting elements come from~$G$-bundles on tori, again reflecting the geometry of~$\Phi_0^{2|1}(X\sq G)$ that is unseen at the level of the de~Rham model. 
\item The moduli stack~$\Phi_0^{2|1}(X)$ possesses a natural family of operators whose determinant constructs the Witten class of~$X$. This gives an analytic construction of the wrong-way maps associated to the string orientation of ${\rm TMF}$ tensored with~$\C$. In this description, the appearance of Eisenstein series in the Witten class comes from traces of powers of the $\bar \partial$-operator on tori. 
\end{enumerate}

A primary goal of this paper is to explain (3) in detail by constructing the Witten class as a holomorphic function on $\Phi_0^{2|1}(X)$, i.e., an element
$$
{\rm Wit}(X)\in \mathcal{O} (\Phi^{2|1}_0(X)):=\mathcal{O} (\Phi_0^{2|1}(X);\omega^{0}),
$$ 
that represents the Witten class under the isomorphism in Theorem~\ref{thethm}.

First we review some facts about the Witten genus. It is a Hirzebruch genus valued in~$\Q((q))$ determined by the characteristic series
\beq
(e^{z/2}-e^{-z/2})\prod_{n\ge 1} \frac{(1-q^n e^{z/2})(1-q^ne^{-z/2})}{(1-q^n)^2}=\exp\left(\sum_{k\ge 1}\frac{E_{2k}(q)}{2k(2\pi i)^{2k}} z^{2k}\right)\label{eq:multWitnonmod}
\eeq
where $E_{2k}$ denotes the (holomorphic) $2k^{\rm th}$ Eisenstein series. The equality above was proved by D.~Zagier~\cite{Zagiermodular}, and reveals the modular properties of the Witten genus conjectured by Witten~\cite{witten_ell} in his original construction. The $2k^{\rm th}$ Eisenstein series is a modular form of weight $2k$ for $k\ge 2$, but~$E_2$ is \emph{not} modular. There is a closely related genus that uses the modular (but non-holomorphic) 2nd Eisenstein series,~$E_{2}^*(q,\bar q)$, with characteristic series
\beq
\exp\left(\frac{E_2^*(q,\bar q)}{2(2\pi i)^2}z^2+\sum_{k\ge 2}\frac{E_{2k}(q)}{2k(2\pi i)^{2k}} z^{2k}\right).\label{eq:multWitnonholo}
\eeq
The associated genus is the \emph{non-holomorphic Witten genus}. For our purposes, we view the genus associated to~(\ref{eq:multWitnonholo}) as taking values in the graded ring of \emph{weak Maass forms}, denoted~${\rm MaF}$, which consists of smooth functions on based lattices that transform as modular forms, but need not be holomorphic. 

\begin{defn} A \emph{(geometric) rational string structure} on an oriented Riemannian manifold $X$ is a 3-form $H\in \Omega^3(X)$ such that $dH=-p_1(X)$ is the first Pontryagin form for the Levi-Civita connection. We call the pair $(X,H)$ a \emph{rational string manifold}. \end{defn}

For rational string manifolds, the characteristic series 
\beq
\sum_{k\ge 2}\frac{E_{2k}(q)}{2k(2\pi i)^{2k}} z^{2k}\label{eq:multWit}
\eeq
defines a Hirzebruch genus with values in the graded ring of weak modular forms, denoted~${\rm MF}$. The characteristic series~(\ref{eq:multWitnonmod}), (\ref{eq:multWitnonholo}) and (\ref{eq:multWit}) define classes in cohomology, called the \emph{non-modular Witten class}, the \emph{non-holomorphic Witten class} and the \emph{Witten class}, respectively denoted $[{\rm Wit}(X)]\in {\rm H}^{\bullet}(X)[[q]]$, $[{\rm Wit}^*(X)]\in {\rm H}^\bullet(X;{\rm MaF})$, and $[{\rm Wit}_H(X)]\in {\rm H}^{\bullet}(X;{\rm MF})$. Although the Witten class does not depend on a choice of rational string structure~$H$ (only its existence) construction of Witten \emph{cocycles} representing these classes will depend on such a choice, so we include it in the notation. The pairing between the fundamental class of $X$ and these Witten classes gives various flavors of the Witten genus of~$X$. 
%

Theorem~\ref{thethm} permits a construction of a cocycle representative for the Witten class for a Riemannian manifold $X$ in terms of a determinant of a family of operators over $\Phi_0^{2|1}(X)$. These operators act on an infinite rank vector bundle~$\mathcal{N}\Phi^{2|1}_0(X)$ over~$\Phi^{2|1}_0(X)$ whose fiber at a map $\phi\colon T^{2|1}\to X$ is the subspace of~$\Gamma(T^{2|1},\phi^*TX)$ that is the orthogonal complement of the constant sections. 
On sections we consider the action of an operator
$$
\Delta^{2|1}_X= \left\{ \begin{array}{cc} -\nabla_{\partial_{\bar z}} \nabla_{\partial_z} + \frac{1}{2}\mathcal{R}(D,D)\nabla_{\partial_z} & {\rm on \ even \ sections} \\ \nabla_{\partial_z} & {\rm on \ odd\ sections} \end{array}\right.
$$
where $\nabla$ denotes the pullback of the Levi-Civita connection on~$X$ along~$\phi$, $\mathcal{R}$ is the pullback of the curvature 2-form of this connection, and $D=\partial_\theta+\theta\partial_{\bar z}$ is an odd vector field on $T^{2|1}$. Note that on even sections, $\Delta^{2|1}_X$ is deformation of the Laplacian on the underlying torus by the curvature of~$X$---a nilpotent deformation. The determinant of this operator involves traces of powers of the curvature operator, which is where the Pontryagin forms of~$X$ enter.

Our first construction of the Witten class is via the $\zeta$-super determinant of the family of operators $\Delta^{2|1}_X$,
$$
{\rm sdet}_\zeta(\Delta^{2|1}(X)):=\frac{{\rm pf}_\zeta((\Delta^{2|1}_X)|_{\rm odd})}{{\rm det}_\zeta((\Delta^{2|1}_X)|_{\rm even})^{1/2}}.
$$ 
It turns out to be convenient to take this super determinant \emph{relative} to the super determinant of an operator denoted $\Delta^{2|1}_n$, whose definition uses the trivial bundle of dimension $n={\rm dim}(X)$ in place of the tangent bundle of~$X$. Along the way to proving Theorem~\ref{thethm} we will show that $C^\infty(\Phi_0^{2|1}(X))= \bigoplus_j \Omega^{j}_{\rm cl}(X;\MaF^{-j})$, so a smooth function on~$\Phi_0^{2|1}(X)$ determines a cocycle in cohomology with coefficients in weak Maass forms. 

A piece of terminology is helpful in the statement below: for a sheaf (or stack) $\mathcal{F}$ on the site of manifolds, a pair of sections $f_0,f_1\in \mathcal{F}(X)$ are \emph{concordant} if there is a section ${f}\in \mathcal{F}(X\times \R)$ such that $i_0^*{f}=f_0$ and $i_1^*{f}=f_1$ where $i_0,i_1\colon X\hookrightarrow X\times \R$ are the inclusions at $0$ and $1$, respectively. By the Yoneda lemma, one can equivalently view a concordance as a smooth homotopy $f\colon X\times \R\to \mathcal{F}$ between the maps of sheaves~$f_0,f_1\colon X\to \mathcal{F}$. 


\begin{thm}\label{thm:Witzeta}
For a choice of Riemannian structure on~$X$, the relative $\zeta$-super determinant of the family of operators $\Delta^{2|1}_X$ gives a function 
$$
\frac{{\rm sdet}_\zeta(\Delta^{2|1}_X)}{{\rm sdet}_\zeta(\Delta^{2|1}_n)}={\rm Wit}^*(X)\in C^\infty({\Phi}_0^{2|1}(X))
$$
that coincides with the non-holomorphic Witten cocycle of~$X$. For any 1-parameter family of Riemannian metrics on $X$, this construction determines a concordance ${\rm Wit}^*(X\times \R)\in C^\infty({\Phi}_0^{2|1}(X\times \R))$ between the representatives of the class $[{\rm Wit}^*(X)]$.

The function ${\rm Wit}^*(X)$ is concordant to the (modular and holomorphic) Witten cocycle if and only if $X$ has a rational string structure, and a choice of rational string structure specifies a concordance~$W_H\in C^\infty(\Phi_0^{2|1}(X\times \R))$ with
$$i_0^*W_H={\rm Wit}^*(X), \quad {\rm and} \quad i_1^*W_H=:{\rm Wit}_H(X)\in \mathcal{O} (\Phi_0^{2|1}(X))\subset C^\infty(\Phi_0^{2|1}(X)).$$
\end{thm}



As a word of warning, there are two independent concordances in the statement of the theorem above: the first comes from changing the metric, and the second comes from a choice of string structure. These both change the cocycle representative of the Witten class, but (typically) in different ways. 

Another approach to the string obstruction and the Witten genus comes from analyzing a Bismut--Freed--Quillen determinant line bundle of~$\Delta^{2|1}_X$ relative to a determinant line of~$\Delta^{2|1}_n$. This ends up being roughly equivalent to the above $\zeta$-determinant story, but unravelling the equivalence clarifies some aspects of these determinant computations over stacks and why (at least in this case) things are unusually simple. This is in contrast, e.g., to the geometrically intricate theory of string structures and obstructions initiated by Freed~\cite{Freed_Det} and Witten~\cite{witten_dirac} in the context of the $2|1$-dimensional sigma model. 

We go into greater detail in~\S\ref{sec:sdet}, but for now we summarize the key steps in generalizing the Bismut--Freed--Quillen construction to the operators $\Delta^{2|1}_X$ over the stack~$\Phi_0^{2|1}(X)$.
One feature simplifies the analysis dramatically: $\Delta^{2|1}_X$ is a family of \emph{invertible} operators, so when evaluated on any $S$-point of~$\Phi_0^{2|1}(X)$ the determinant line bundle on~$S$ of the pullback of~$\Delta^{2|1}_X$ is canonically trivialized by its determinant section. However, this line bundle is not canonically trivialized metrically; the norm squared of the trivializing section $\sigma$ is the pullback of 
\beq
\|\sigma\|^2=\frac{{\rm sdet}_\zeta((\Delta^{2|1}_X)^*\Delta^{2|1}_X)}{{\rm sdet}_\zeta((\Delta^{2|1}_n)^*\Delta^{2|1}_n)}\label{eq:normsigma}
\eeq
to~$S$.
One can ask whether there is a choice of metric trivializations at each $S$-point that descends to a trivialization over the stack, i.e., a square root of $\|\sigma\|^2$ natural in~$S$ that rescales the line over each~$S$ so that the trivializing section has unit norm. Naturality makes this equivalent to asking whether the function $\|\sigma\|^2$ on $\Phi_0^{2|1}(X)$ has a square root, which a priori need not exist nor be unique. Conveniently, the $\zeta$-determinant computation underlying Theorem~\ref{thm:Witzeta} gives an essentially unique square root.

\begin{lem} ${\rm Wit}^*(X)$ is the unique function on $\Phi_0^{2|1}(X)$ such that
\begin{enumerate}
\item the norm of ${\rm Wit}^*(X)$ satisfies (\ref{eq:normsigma});
\item at any $S$-point of $\Phi_0^{2|1}(X)$, the restriction of ${\rm Wit}^*(X)$ to the reduced manifold of $S$ is the function~1. 
\end{enumerate}\label{lem:uniquefun}
\end{lem}

\begin{rmk} There are different universal choices of square root of $\|\sigma\|^2$ that give rise to \emph{non-isomorphic} (and possibly nontrivial) determinant line bundles over the stack~$\Phi_0^{2|1}(X)$. The issue is that although the norm of~$\sigma$ is fixed, its phase is not. This phase can change over isomorphisms between $S$-points, and descending a line bundle over $S$-points of~$\Phi_0^{2|1}(X)$ to the stack itself requires one choose a regularization procedure for these phases. The lemma above picks out one such regularization, but there are others. 
\end{rmk}

In choosing the metric trivialization that rescales by ${\rm Wit}^*(X)^{-1}$ at each $S$-point, we identify the super determinant section $\sigma$ with the function~${\rm Wit}^*(X)$ on~$\Phi_0^{2|1}(X)$. 
A computation finds~${\rm Wit}^*(X)$ is not invariant under the required infinitesimal anti-holomorphic dilations of tori to define a cocycle, which boils down to the fact that $E_2^*$ isn't a holomorphic function on lattices. However, at each $S$-point of $\Phi_0^{2|1}(X)$ there is a canonical modification to the phase of the determinant to make the pullback of~${\rm Wit}^*(X)$ holomorphic, namely
\beq
\exp\left(\frac{-p_1(TX)}{2\pi \cdot \vol}\right)\cdot {\rm Wit}^*(X)={\rm Wit}(X),\label{eq:p1mod}
\eeq
where $\vol$ is the function on lattices that reads of the volumes of tori. This constructs the non-modular Witten class. In terms of the geometry of $\Phi_0^{2|1}(X)$, the failure of modularity means that the functions at each $S$-point determined by~${\rm Wit}(X)$ do not descend to a function on~$\Phi^{2|1}_0(X)$, but instead assemble into a section of a line bundle~$\mathcal Str$. Indeed, this modified class ${\rm Wit}(X)$ is a cocycle representative of the twisted fundamental class of~$X$ in~$\TMF(X)\otimes \C$, where the twist comes from the degree~4 cohomology class $p_1(TX)$. A choice of (rational) string structure is equivalent to a null-concordance of~$\mathcal Str$.

\begin{thm} The line bundle $\mathcal{S}tr$ is concordant to the trivial line bundle if and only if $X$ has a rational string structure. A choice of rational string structure specifies such a concordance, namely a line bundle~$\widetilde{\mathcal{S}tr}$ on $\Phi_0^{2|1}(X\times \R)$ with a section $\tilde{\sigma}$ whose restriction to the source of the concordance is~${\rm Wit}(X)$ and the restriction to the target is the function
$$
i_1^*\tilde{\sigma}={\rm Wit}_H(X)\in \mathcal{O} (\Phi_0^{2|1}(X)).
$$ 

\label{thm2}
\end{thm}

\subsection{Physical motivation from 1-loop quantization}\label{sec:introphys}

The constructions underlying Theorems~\ref{thethm} and~\ref{thm2} can be motivated by standard maneuvers in physics applied to the $2|1$-dimensional sigma model with target a Riemannian manifold~$X$. We give a brief overview here, with more details in \S\ref{sec:21motivate} and \S\ref{sec:21Hess}; see also~\S\ref{sec:SUSYmechinter} and \S\ref{sec:11Hess} for the $1|1$-dimensional case. Treating the space of fields as a \emph{stack} throughout the analysis is a relatively new (though widely-anticipated) ingredient in the story. This refinement is crucial for our cocycle model---see Remark~\ref{rmk:stack} for a toy example highlighting the distinction. 

We start with a general recipe from physics called \emph{1-loop quantization}; e.g., see~\cite{Alvarez} and \cite{WitteninStrings} in the case of $1|1$-dimensional supersymmetric quantum mechanics. Given a space (or stack) of fields $\mathcal{F}$ with a classical action functional~$\mathcal{S}$, the subspace on which~$d\mathcal{S}=0$ are called \emph{classical solutions}, denoted~$\mathcal{F}_{\rm cl}\subset\mathcal{F}$. There is a function on these classical solutions called the~\emph{energy}, and the subspace on which energy is minimized consists of the \emph{classical vacua}, denoted $\mathcal{F}_{\rm vac}\subset\mathcal{F}_{\rm cl}$. A crucial point is that classical vacua usually form a finite-dimensional subspace of~$\mathcal{F}$, making the geometry of~$\mathcal{F}_{\rm vac}$ easier to study. For example, standard differential geometry can be used to construct volume forms on~$\mathcal{F}_{\rm vac}$, whereas this is notoriously difficult on the space~$\mathcal{F}$ which is typically infinite-dimensional. With suitable non-degeneracy assumptions on~$\mathcal{S}$, the Hessian of the classical action restricted to vacua determines a family of invertible operators over~$\mathcal{F}_{\rm vac}$ called \emph{kinetic operators} that act on sections of the normal bundle to the inclusion~$\mathcal{F}_{\rm vac}\subset\mathcal{F}$. In this situation, the determinant of the kinetic operators yields a non-vanishing function on classical vacua. A volume form on~$\mathcal{F}_{\rm vac}$ (which is often canonical) can by modified by this function. The resulting total volume of~$\mathcal{F}_{\rm vac}$ is the \emph{1-loop partition function}, as it computes the contribution to the partition function from Feynman diagrams with 1-loop, e.g., compare Etinghof \cite{EtingofQFT} Theorem~3.5. 
For a class of supersymmetric field theories, this 1-loop approximation to the partition function is expected to coincide with a hypothetical non-perturbative partition function by a formal application of the Atiyah--Bott fixed point theorem~\cite{witten_dirac,WitteninStrings}. 

Specializing to our main example, the $2|1$-dimensional sigma model has as fields maps $\phi\colon T^{2|1}_\Lambda\to X$ from super tori to~$X$. The classical action is invariant under a certain isometry group of~$T^{2|1}$, so can be promoted to a function on a stack of fields, denoted~$\Phi^{2|1}(X)$, whose objects are maps~$\phi$ and morphisms are commuting triangles
\beq
\begin{tikzpicture}[baseline=(basepoint)];
\node (A) at (0,0) {$T_\Lambda^{2|1}$};
\node (B) at (3,0) {$T_{\Lambda'}^{2|1}$};
\node (C) at (1.5,-1.5) {$X$};
\node (D) at (1.5,-.6) {$\circlearrowright$};
\draw[->] (A) to node [above=1pt] {$\cong$} (B);
\draw[->] (A) to node [left=1pt]{$\phi$} (C);
\draw[->] (B) to node [right=1pt]{$\phi'$} (C);
\path (0,-.75) coordinate (basepoint);
\end{tikzpicture}\nonumber
\eeq
where the horizontal arrow is an isometry between super tori. Critical points of the action define a substack of classical solutions, $\Phi^{2|1}_{\cl}(X)\subset \Phi^{2|1}(X)$. The energy function on classical solutions in this case is the standard energy of the underlying map from an ordinary torus to~$X$; in particular, it is nonnegative. The stack $\Phi_0^{2|1}(X)\subset \Phi^{2|1}_{\cl}(X)$ is the stack of \emph{classical vacua}, consisting of the classical solutions with zero energy; see~\S\ref{sec:21motivate}. One important feature unique to the $2|1$-dimensional sigma model is that the stack~$\Phi_0^{2|1}(X)$ is in a sense topological: it is independent of the metric on~$X$ chosen to define the classical action. 

The restriction of the classical action vanishes to second order on $\Phi_0^{2|1}(X)$ with its Hessian being encoded by the \emph{kinetic operators}~$\Delta^{2|1}_X$ (see~\S\ref{sec:21Hess})
$$
\Hess(\mathcal{S})(\nu,\nu')=\langle \nu ,\Delta^{2|1}_X\nu'\rangle,
$$
for sections $\nu,\nu'$ of the normal bundle to the inclusion $\Phi_0^{2|1}(X)\subset \Phi^{2|1}(X)$. Then 1-loop quantization takes a regularized determinant of $\Delta^{2|1}_X$. However, in infinite dimensions determinants are not functions, but rather sections of a line. So while this quantization is unobstructed locally (as lines can always be trivialized), there can be non-trivial effects globally over~$\Phi_0^{2|1}(X)$. Relatedly, the line bundle might be trivializable, but a trivializing global section may not possess certain desired properties, e.g., holomorphicity. Often there is a modification to this section that is local in~$X$ and restores some of these desired properties; the relevant modification in this paper is~(\ref{eq:p1mod}). However, in this case the result is no longer a function but instead a section of a line bundle. Without some choices of additional structure on~$X$ whose existence can be obstructed topologically---i.e., a rational string structure---there is no canonical trivialization of this line.

These types of global effects with topological obstructions are examples of \emph{anomalies}, and have been mathematically formalized in the language developed by Bismut and Freed~\cite{BismutFreed1,Freed_Det} in the context of the (non-perturbative) $2|1$-dimensional sigma model. In this case, the anomaly is a line bundle~$\mathcal{L}$ with metric and connection over the double loop space of~$X$, that is defined as the Pfaffian of a family of twisted Dirac operators~\cite{Freed_Det,BunkeString}. The curvature 2-form of~$\mathcal{L}$ is $\pi_!\circ \ev^*p_1(TX)/2$, where $\ev^*$ is the pullback along evaluation and $\pi_!$ is the pushforward (i.e., integration) along the projection,
\beq
{\sf Map}(\T^2,X)\stackrel{\pi}{\leftarrow} \T^2\times{\sf Map}(\T^2,X)\stackrel{\ev}{\to} X.\label{eq:sigmaev}
\eeq
Physically~$\mathcal{L}$ and its Pfaffian section result from the path integral over fermonic fields in the sigma model; since the result is a section of a line bundle, one can't even formulate the \emph{integrand} for the remaining integral over the bosonic fields~${\sf Map}(\T^2,X)$, much less make sense out of the integral. Metric trivializations of~$\mathcal{L}$ allow one to view the Pfaffian section as a function, which then (at least formally) one can integrate. If we restrict (\ref{eq:sigmaev}) to the constant maps, the class $\pi_!\circ \ev^*p_1(TX)$ is identically zero, and the associated line bundle is \emph{canonically} trivialized. The perturbative approximation studies a small neighborhood of these constant maps, so we should expect this line bundle to trivialize topologically; however, it can still have nonvanishing curvature. Another obstruction to triviality comes from automorphisms of fields, and one really wants an equivariant trivialization of~$\mathcal{L}$  over the double loop \emph{stack}. Although similar obstructions to trivializations of~$\mathcal Str$ show up in this paper, we don't know of a direct connection between Bismut and Freed's (non-perturbative) anomaly theory and the perturbative one. 

In addition to this 1-loop quantization story, the stack $\Phi_0^{2|1}(X)$ can be viewed as a subcategory of a bordism category over~$X$; see~\S\ref{sec:Stolz--Teichner}. From this point of view, families of field theories over~$X$ can be evaluated to produce (partition) functions on~$\Phi_0^{2|1}(X)$. Hence,~$\Phi_0^{2|1}(X)$ plays a double role: as the classical vacua for a particular theory, and as the world sheet for families of (possibly quantum) field theories. Functions arising from families of field theories satisfy an additional condition, namely supersymmetric partition functions depend \emph{holomorphically} on tori. This comes from the fact that the value on a torus is itself a trace of a linear map associated to a cylinder, and the super geometry of cylinders and their gluings reveals that the trace depends holomorphically on the modulus of the cylinder. Theorem~1.15 of~\cite{ST11} makes this precise for~$X=\pt$ using supersymmetric cancellation arguments analogous to those in physics proofs of the Atiyah--Singer index theorem. This holomorphic dependence is the distinction between $\Gamma(\Phi_0^{2|1}(X);\omega^{\bullet/2})$ and $\mathcal{O}(\Phi_0^{2|1}(X);\omega^{\bullet/2})$, i.e., between all sections and those that define cocycles. 

\subsection{Comparison with Witten's Dirac operator on loop space}\label{sec:introWit}

In his original construction, Witten formally applies the Atiyah--Bott fixed point theorem~\cite{AtiyahBott1} to compute the index of a would-be Dirac operator on the free loop space of~$X$ with an $S^1$-action from rotating loops~\cite{witten_dirac}. Formally, this gives a Dirac operator twisted by an infinite-dimensional vector bundle, constructed out of the normal bundle to the inclusion~$X\subset LX$ of the constant loops, or equivalently, the $S^1$-fixed points. The fibers of this bundle carry an $S^1$-action with finite-dimensional weight spaces, so Witten repackages his Dirac operator as one twisted by an infinite sequence of finite-dimensional vector bundles. Taking the index term-by-term in the sequence gives an integral power series in a variable~$q$, with the power of $q$ encoding the $S^1$-action. This is the Witten genus of~$X$. Applying physical intuition, Witten claims that this genus is the $q$-expansion of a modular form provided that~$X$ has a string structure, and this was verified by Zagier~\cite{Zagiermodular}. Bott and Taubes' clarified the geometry of Witten's construction in terms of fixed point formulae along the way to proving the rigidity of the Witten genus~\cite{BottTaubes}. 

In physical terms, Witten's construction is the quantization in the Hamiltonian formalism of the same classical field theory as ours. Morally this requires that one study all maps of $2|1$-dimensional super cylinders to~$X$, form an infinite-dimensional symplectic super manifold of classical solutions, and apply methods of geometric quantization. However, this classical field theory has a type of perturbative approximation whose fields are constant maps~$x$ of super cylinders to~$X$ together with a section of of the pullback bundle,~$x^*TX$. This approximate space of fields can also be described as the normal bundle to fields invariant under the $S^1$-action that rotates super cylinders, which is the point of view that leads to Witten's formal application of the Atiyah--Bott fixed point theorem. Geometric quantization of sections of $x^*TX$ produces a family of Hilbert spaces and (super) time-evolution operator parametrized by~$X$. The fiber at a fixed~$x$ is Witten's infinite sequence of twisted spinor bundles, and time-evolution is encoded by the $S^1$-action together with the sequence of twisted Dirac operators. The reason why this type of quantization is relatively straightforward mathematically is that for a fixed constant map~$x$ the classical field theory in question is free: fields are maps into the tangent space, and the relevant action functional is purely quadratic. Geometric quantization of free theories is essentially algebraic, being a Fock space construction depending on a choice of polarization. Although this is unobstructed locally in~$X$, picking polarizations in families can present difficulties, which is exactly where the spin structure on~$X$ comes in for Witten's construction. 

In this paper we apply a similar perturbative approximation to the $2|1$-dimensional sigma model, but in the Lagrangian formalism. Indeed, a different way of describing 1-loop quantization is that the Hessian of the classical action defines a family of free field theories over classical vacua. The fields are sections of the normal bundle to the $\T^2$-fixed point set for the translational action of $\T^2$ on super tori, and the classical action is determined by the kinetic operator. At each point, quantization is unobstructed: in this easy case, quantization amounts to taking the determinant of an operator. However, just as in the Hamiltonian case, obstructions to quantization can arise in families via determinant line bundles, as discussed in the previous subsection. 

The Hamiltonian approach taken by Witten makes clear the integrality properties of his genus, but (without appealing to physical intuition) modularity seems miraculous and must be verified computationally after the fact~\cite{Zagiermodular}. On the other hand, although the Lagrangian approach taken below doesn't detect integrality properties, it constructs the Witten genus over the moduli stack of tori. This provides some geometric and physical intuition for modularity of the Witten genus of a string manifold. 

In the case of $1|1$-dimensional supersymmetric quantum mechanics, the relation between 1-loop quantization and Hamiltonian quantization is well-known, being Alvarez-Gaum\'e's physical proof of the index theorem~\cite{Alvarez}. We also drew inspiration from a more recent accounting by Witten~\cite{WitteninStrings}. Below we repackage these sorts of physical computations in a different conceptual framework (e.g., that makes the connection to differential cocycles more apparent), and we generalize them to the $2|1$-dimensional case. As discussed in \S\ref{sec:introphys}, this generalization one dimension up confronts some new features related to anomalies and the rational string obstruction. This integrates some subtle aspects of the physical story with expected features on the topological side.

\subsection{Super-generalizations of geometric index bundles and determinant lines} \label{sec:reldetline}

The main objects of study in this paper can be understood as super-geometric generalizations of standard objects in index theory. The stack of classical vacua,~$\Phi_0^{2|1}(X)$, is closely related to the kernel bundle associated to a geometric family of spin tori in the sense of Bismut--Freed~\cite{BismutFreed1,BismutFreed2} equipped with an action induced by spin isometries between tori. The stack $\Phi_0^{2|1}(X)$ of classical vacua generalizes this standard object by considering the action by \emph{super Euclidean} isometries; see Example~\ref{ex:superfication} for a comparison between spin and super Euclidean geometries. The other main player, the family of operators~$\Delta^{2|1}_X$, are deformations of (Dirac) Laplacians on this same geometric family of spin manifolds. The new ingredient here is a family of nilpotent deformations of these Laplacians coming from odd directions in the stack~$\Phi_0^{2|1}(X)$. 

In more detail, consider the manifold~$L$ of based, oriented lattices in~$\C$ (note $L\cong \C^\times\times\mathfrak{h}$). One presentation of the moduli stack of elliptic curves is the quotient stack~$L\sq \C^\times\times \SL_2(\Z)$, and pulling back the universal curve along the projection $L\to L\sq \C^\times\times\SL_2(\Z)$ gives a natural family of geometric spin manifolds over~$L$, where we endow the elliptic curves with the standard flat metric and odd (or nonbounding) spin structure. We can further pull this family back along the projection~$L\times X\to L$, resulting in a geometric family of spin tori $F\to L\times X$. There is a vector bundle over~$F$ gotten from pulling back the super vector bundle~$\Pi TX\to X$ along the projection $F\to X$. The usual Bismut--Freed setup leads to a family of twisted Dirac operators parametrized by~$L\times X$; in each fiber, this is the Dirac operator on a torus twisted by the fiber of the tangent bundle of~$X$, viewed as a purely odd vector space. The kernel bundle is isomorphic to $L\times \Pi TX\to L\times X$, since in this case harmonic spinors are just constant sections. Note, however, that this identification chooses a trivialization of the spinor bundle on the tori. As we shall see, the stack $\Phi_0^{2|1}(X)$ has a groupoid presentation with objects~$L\times \Pi TX$. 

Naturality of the kernel bundle with respect to isometries of geometric families of spin manifolds gives rise to an action by an isometry group $\SL_2(\Z)\times \T^2_L\ltimes \C^\times$, where $\SL_2(\Z)$ changes the lattice, $\T^2_L$ acts by translations, and $\C^\times$ acts by dilations on tori. Note that $\SL_2(\Z)\times \T^2_L$ acts on tori by isometries preserving the metric, whereas $\C^\times$ does not. Although this dilation action does not preserve the geometric structures, the kernel bundle is modified in a predictable way (namely, dilation on the fibers of~$\Pi TX$), which allows one to view the total space of the kernel bundle as a stack whose objects are $L\times \Pi TX$ and whose morphisms are $\SL_2(\Z)\times \T^2_L\ltimes \C^\times$. The stack~$\Phi_0^{2|1}(X)$ extends this construction by enlarging the space of maps between tori to a certain \emph{super} isometry group, which roughly adds an odd square root of the anti-holomorphic translations in~$\T^2_L$. 

The operators $\Delta_X^{2|1}$ come from a family of deformations of the Laplacian on the family~$F$, with deformations parametrized by the fibers of the kernel bundle~$L\times \Pi TX\to L\times X$. For the original family, we take Laplacians on the family of tori twisted by the pullback of~$TX$, and we restrict to the orthogonal complement of their kernel, i.e., we start with a family of invertible operators on~$L\times X$. The deformation~$\Delta^{2|1}_X$ can be motivated physically, but basically uses the fact that the curvature operator is an order zero differential operator over~$\Pi TX$: it is a 2-form valued in endomorphisms of~$TX$, which can be repackaged as a quadratic function on the supermanifold~$\Pi TX$ valued in endomorphisms. Furthermore, this operator is nilpotent (since its differential form degree is positive), so $\Delta^{2|1}_X$ is a nilpotent deformation of a family of invertible operators, and hence is a family of invertible operators. This leads to a canonically trivialized Bismut--Freed--Quillen (super) determinant line bundle on~$L\times \Pi TX$. Understanding equivariant trivializations of this line bundle over $\Phi_0^{2|1}(X)$ that are appropriately holomorphic leads to an analysis of rational string structures on~$X$. 

\subsection{Relation to the Segal--Stolz--Teichner program}\label{sec:Stolz--Teichner}

G.~Segal proposed a connection between 2-dimensional field theories and elliptic cohomology~\cite{segal_elliptic}. The fundamental objects of study are field theories over a manifold $X$. These are linear representations of a category whose morphisms are bordisms equipped with a smooth map to~$X$. With supersymmetry being such a crucial ingredient in Witten's construction of his genus, Stolz and Teichner constructed a super Euclidean bordism category over a smooth manifold~$X$, denoted $2|1\EB(X)$~\cite{ST04,ST11}. The objects consist of closed $1|1$-dimensional super manifolds with a map to~$X$, and the morphisms are compact $2|1$-dimensional super manifolds with super Euclidean structure and a map to~$X$. To account for the symmetries super Euclidean manifolds possess, this bordism category is defined as being internal to stacks on the site of super manifolds, i.e., the collections of objects and morphisms each form stacks. As usual, disjoint union of bordisms gives this category a symmetric monoidal structure. Stolz and Teichner further define a symmetric monoidal functor ${\sf T}$, called a \emph{twist}, from this bordism category to the category ${\sf Alg}$; this target has a stack of objects consisting of topological algebras and algebra automorphisms and a stack of morphisms consisting of bimodules and bimodule automorphisms. They consider natural transformations~$F$
$$
\begin{tikzpicture}[node distance=2cm, auto]
  \node at (0,0) (A) {$2|1\EB(X)$};
  \node at (2.15,0) (B)  {$\Downarrow$ F};
  \node at (4,0) (C)  {${\sf Alg},$};

\draw[->,bend left] (A.20) to node {${\sf T}^{\otimes k}$} (C.160);
\draw[->, bend right] (A.340) to node [swap] {$\one$} (C.200);
\end{tikzpicture}\label{eq:twist}
$$
denoting the category of such natural transformations by $2|1\EFT^k(X)$ for $k\in \Z$. Roughly, these are field theories valued in modules over algebras~$A_k$,~$k\in \Z$ such that $A_k\otimes A_l\cong A_{k+l}$. 

An important observation of Stolz and Teichner's~\cite{ST04} is that for such field theories to have a chance of admitting Mayer--Vietoris sequences, they must be fully local in the higher-categorical sense. The essential idea is that a map of a circle---viewed as an object of a 2-dimensional bordism category---into a smooth manifold~$X=U\bigcup V$ will generally not map to one of $U$ or $V$. Mayer--Vietoris sequences would require field theories over $U$ and $V$ to determine a field theory over $X$, which effectively forces one to associate data to a chopped up circle, i.e., to intervals. In short, we need to require field theories over~$X$ to satisfy a higher sheaf condition; in a modern language, 2-dimensional field theories need to form a 2-stack on~$X$.  Stolz and Teichner conjecture the existence of a higher-categorical refinement of their definition of $2|1\EFT^k(X)$ such that there is a natural ring isomorphism
$$
{\rm TMF}^k(X)\cong 2|1\EFT^k(X)/{\rm concordance} \quad\quad\quad {(\rm conjectural}).
$$

However, there is cheaper way to enforce locality: we can require the maps from bordisms to~$X$ to be \emph{constant}, so that such bordisms over $X=U\bigcup V$ automatically lie in either $U$ or~$V$. When the bordisms in question are tori, a modification of this idea incorporating the relevant super geometry leads to our stack of classical vacua,~$\Phi^{2|1}_0(X)$. In this way, one can view our cocycles for ${\rm TMF}\otimes \C$ as extracting a sheaf of sets (or \emph{0-categories}) from the Stolz--Teichner presheaf of 1-categories $X\mapsto 2|1\EFT(X)$.

Any candidate definition of the extended $2|1$-dimensional bordism category will be far more intricate that~$\Phi_0^{2|1}(X)$, but analyzing the restriction of these would-be extended field theories to the subcategory~$\Phi_0^{2|1}(X)$ gives a glimpse at what information such a higher-categorical object might contain. Restriction of a degree~$k$ twisted field theory to the subcategory of $2|1\EB_{\rm c}(X)\subset 2|1\EB(X)$ of closed, connected tori equipped with a map to~$X$ picks out a section of a line bundle~$\mathcal{L}_{\sf T}^{k}$ over $2|1\EB_{\rm c}(X)$, where $\sf T$ determines~$\mathcal{L}_{\sf T}$. When $X=\pt$, Stolz and Teichner constructed a map $Z$ (as a consequence of Theorem 1.15 in~\cite{ST11}),
\beq
\begin{array}{lllllllll}
 2|1\EFT^k(X)&\to& \Gamma(2|1\EB_{\rm c}(X); \mathcal{L}_{\sf T}^{ k})&\stackrel{Z}{\to}&\Gamma  (\Phi_0^{2|1}(X);\omega^{k/2}),
 \end{array}\label{eq:CH}
\eeq
and for general~$X$ this is constructed in~\cite{DBEEll}. In this way, Theorem \ref{thethm} gives a map from $2|1\EFT^k(X)$ to cocycles for ${\rm TMF}(X)\otimes \C$ independent of the choice of higher-categorical refinement under investigation by Stolz and Teichner. This gives a link between $2|1$-Euclidean field theories and elliptic cohomology independent of these subtleties. 

Analyzing the image of~(\ref{eq:CH}) illuminates the manner in which field theories refine cocycles for ${\rm TMF}(X)\otimes \C$. Some of the 1-categorical aspects of this refinement have been worked out using a category of energy zero super circles and super annuli over~$X$, denoted~${\sf Ann}^{2|1}_0(X)$, which was studied for $X=\pt$ in~\cite{PokmanPhD,ST11} and over general~$X$ in~\cite{DBEEll}. Representations of this category refine functions on~$\Phi_0^{2|1}(X)$ in much the same way that a group representation refines its character. Among all functions on $\Phi_0^{2|1}(X)$, those that come from representations of ${\sf Ann}^{2|1}_0(X)$ satisfy an integrality condition, which when $X=\pt$ is exactly the refinement from a modular form to an \emph{integral} modular form. This follows from identifying the coefficients in the $q$-expansion of a modular form as dimensions of vector spaces. Over general~$X$, the integrality is akin to the way that a vector bundle with connection refines a its Chern form. An open question in the Stolz--Teichner program is to find an enhancement of these field theories in which the existence and uniqueness of 2-categorical refinements of cocycles mimics the the complicated torsion present in topological modular forms. 

We now turn to how our construction of the Witten genus fits into this story. In short, it can be viewed as a piece of a conjectured Riemann--Roch or local index theorem for TMF. The local index theorem in K-theory is summarized by the commutative diagram
\begin{equation}
\begin{array}{c}
\begin{tikzpicture}[node distance=3.5cm,auto]
  \node (A) {${\rm K}^{k+n}(X)$};
  \node (B) [node distance= 4.5cm, right of=A] {$ {\rm K}^{k+n}(X)\otimes \C$};
  \node (C) [node distance = 1.5cm, below of=A] {${\rm K}^{k}(\pt)$};
  \node (D) [node distance = 1.5cm, below of=B] {${\rm K}^k(\pt)\otimes \C$};
  \draw[->] (A) to node {${\rm ch}$} (B);
  \draw[->] (A) to node [swap] {$\pi_!^{\rm top},\ \pi_!^{\rm an}$} (C);
  \draw[->] (C) to node {${\rm ch}$} (D);
  \draw[->] (B) to node {$\pi_!(-\smallsmile [\hat{A}(X)]), \ \int_X(-\wedge \hat{A}(X))$} (D);
\end{tikzpicture}\end{array}
\label{diag:red}
\end{equation}
where $\hat{A}(X)$ denotes the $\hat{A}$-form associated to a choice of metric and Levi-Civita connection on $X$, and for each  downward arrow there are two maps that are equal, namely the topological pushforward and the analytic pushforward. We have used ${\rm K}(X)\otimes \C\cong {\rm HP}_{\rm dR}(X)$ where ${\rm HP}_{\rm dR}$ denotes 2-periodic de~Rham cohomology. 

Stolz and Teichner have described an analogous picture that generalizes the left side of~(\ref{diag:red}) to~TMF. They conjecture the existence of a quantization functor for $2|1$-EFTs, 
$$
\pi_!^{\rm an}\colon 2|1\EFT^k(X)\to 2|1\EFT^{k-n}(\pt)\quad\quad {\rm (conjectural)}
$$ 
that gives an analytical index that is equal to the topological index, $\pi_!^{\rm top}$, defined by the Ando, Hopkins, Strickland and Rezk \cite{AHS,AHR}: evaluating $\pi_!^{\rm an}$, on concordance classes is conjectured to give a map equal to~$\pi_!^{\rm top}$ evaluated on cohomology classes. This situation is summarized by the left face of the cube in Figure~\ref{fig1}. Using the generalized Chern character for~TMF, one can also incorporate a local index. Theorem~\ref{thm2} is an analytic construction that agrees with the topological pushforward in~${\rm TMF}\otimes \C$, fitting into the right face of the cube in Figure~\ref{fig1}. More importantly, the analytic techniques of the construction come from a quantization procedure for~$2|1$-dimensional field theories, lending credence to Stolz and Teichner's much deeper conjecture. 

\begin{figure}
$$
\begin{tikzpicture}[node distance=4cm, auto,descr/.style={fill=white}]
[
  \node (A) {$2|1\EFT^{k+n}(X)$};
  \node [right of=A,node distance=5.4cm] (B) {$\mathcal{O}  (\Phi_0^{2|1}(X);\omega^{\bullet/2})$};
\node [below of=A] (C) {$2|1\EFT^{\bullet-n}(\pt)$};
\node [below of=B] (D) {$\mathcal{O}  (\Phi_0^{2|1}(\pt);\omega^{\bullet-n/2})$};

  \node (A1) [right of=A, above of=A, node distance=3.1cm] {${\rm TMF}^{\bullet}(X)$};
  \node [right of=A1] (B1) {${\rm TMF}^{\bullet}(X)\otimes \C$};
    \node [below of=A1] (C1)  {${\rm TMF}^{\bullet-n}(\pt)$};
\node [below of=B1] (D1) {${\rm TMF}^{\bullet-n}(\pt)\otimes \C$};

 \draw[->] (B1) to node {$\pi_!(-\smallsmile[{\rm Wit}(X)])$} (D1);
 \draw[->,dashed] (A) to (A1);
 \draw[->,dashed] (A) to node [swap] {$\pi_!^{\rm an}$}  (C); 
 \draw[->] (A1) to node [swap] {$\pi_!^{\rm top}$} (C1);
 \draw[->] (A1) to node {${\rm ch}$} (B1);
 \draw[->] (D) -- (D1);
 \draw[double distance=4pt, white] (B) -- (B1);
  \draw[->] (B) -- (B1);
  \draw[->] (C1) to node {${\rm ch}$} (D1);
 \draw[double distance=4pt,white] (B) to node {$\phantom{\int(-\wedge{\rm Wit}(X))}$} (D);
   \draw[->] (B) to node[descr] {$\int(-\wedge{\rm Wit}(X))$} (D);
  \draw[->] (C) to node {$\partition$} (D);
 \draw[->,dashed] (C) to (C1); 
  \draw[double distance=4pt,white] (A) -- (B);
 \draw[double distance=4pt,white] (A) to node {$\phantom{\partition}$} (B);
 \draw[->] (A) to node {$\partition$} (B);

\end{tikzpicture}
$$
\caption{The above diagram summarizes a conjectured local index theorem for TMF using supersymmetric field theories. The solid arrows have been constructed, whereas the dotted arrows remain conjectural. The front face describes a local analytical index whereas the back face describes a local topological index so that commutativity of the cube would give a TMF-analog of the local Atiyah--Singer index theorem. The arrows from the front face to the back face arise from taking concordance classes. }\label{fig1}
\end{figure}


\subsection{Relation to Costello's construction of the Witten genus}\label{intro:cost}

K. Costello constructs the Witten genus~\cite{costello_WG1,costello_WG2} using a quantum field theory based on a (formal, derived) stack of maps from an elliptic curve to a complex manifold~$X$ that are nearly constant in the language of Gelfand--Kazhdan formal geometry. Costello then interprets the Witten class as a projective volume form on this derived stack, so it defines a type of integration (or pushforward) for functions on the stack. Below, we consider maps from super tori to a smooth manifold~$X$ that are nearly constant in a super-geometric sense, and our construction of the Witten class determines a pushforward on functions. Very roughly, the pair of classical field theories that construct the Witten class are different ways of encoding the same problem in deformation theory, namely deformations of a constant map from an elliptic curve to~$X$ with a constant section of the pullback tangent bundle into a non-constant map from an elliptic curve to~$X$ and a non-constant section of the pullback tangent bundle. There is a duality between formal deformation problems and differential graded Lie algebras, as developed by Deligne, Drinfeld and Feigin and thereafter by Kontsevich--Soibelman, Lurie, Manetti and others. In our setting, the sheaf of sections of the normal bundle for the inclusion~$\Phi_0^{2|1}(X)\subset \Phi^{2|1}(X)$ defines a sheaf of deformations: a section of the normal bundle deforms a map using the Riemannian exponential map on~$X$. Associated to this sheaf of deformations is a sheaf of differential graded Lie algebras that is close in spirit to the sheaf of $L_\infty$ algebras studied by Costello.


In addition to these structural similarities, our construction shares many computational features with Costello's. The determinant of the Hessian of the action functional is exactly the 1-loop contribution to the partition function computed using Feynman diagrams. Costello shows that quantizing his theory while preserving certain symmetries is such a 1-loop calculation. The precise form of Costello's action functional is quite different from ours, though one might hope for a relationship in a large-volume limit of our theory. 

An important difference between the approaches is that Costello requires the target manifold be K\"ahler, whereas our construction applies to all oriented smooth manifolds. There seems to be an underlying conceptual reason for this, coming from the geometry of 2-dimensional sigma models (for example, see~\cite{Witten0_2}): the $2|1$-dimensional sigma model with target a K\"ahler manifold $X$ has extra symmetry, owing to the fact that the single supersymmetry breaks into two under the holomorphic and anti-holomorphic decomposition of $T_\C X$. This allows one to perform a \emph{half-twist}\footnote{There is an unfortunate collision of terminology: this notion of twist has nothing to do with the Stolz--Teichner twists described in the previous subsection.} of the theory which produces a square zero odd symmetry, leading to constructions in derived geometry (in particular, AKSZ models) utilized by Costello. Physical arguments lead one to expect a special subalgebra of the algebra of observables---called \emph{chiral differential operators} by V.~Gorbounov, F.~Malikov and V.~Schechtman \cite{GMS}---that are conformally invariant. Indeed, Costello expects the factorization algebra produced by his quantization procedure to be equivalent to these chiral differential operators. In the case of a general target manifold we expect the partition function to be conformally invariant, but there need not be a nontrivial subalgebra of conformally invariant observables. Hence, the factorization algebra of observables ought to be more complicated when the target manifold is not K\"ahler. 

\subsection{Notation, terminology and some background}\label{terminology}

Throughout, $X$ will denote an oriented, closed, smooth manifold, and ${\sf Mfld}$ will denote the category of smooth manifolds and smooth maps. For a smooth manifold~$X$, $C^\infty(X)=C^\infty(X,\C)$ will denote the algebra of complex-valued smooth functions on~$X$. 

\subsubsection{Super manifolds} We take a \emph{$k|l$-dimensional super manifold} to be a locally ringed space whose structure sheaf is locally isomorphic to $C^\infty(U)\otimes_\C\Lambda^\bullet(\C^l)$ as a super algebra over~$\C$ for $U\subset \R^k$ an open submanifold. These are sometimes called $cs$-manifolds (e.g.,~\cite{DM}). We denote the category of super manifolds and maps of super manifolds by~${\sf SMfld}$. For any super manifold~$N$, there is a \emph{reduced manifold} we denote by $|N|$ and a canonical map $|N|\hookrightarrow N$ induced by the map of superalgebras $C^\infty(N)\to C^\infty(N)/I\cong C^\infty(|N|)$ where $I$ denotes the ideal of nilpotent elements in the structure sheaf of $N$. We will use notation like $z,\bar{z}$ or $f,\bar{f}$ for elements of $C^\infty(N)$ that are complex conjugates in their image under the quotient $C^\infty(N)\to C^\infty(|N|)$. By M.~Batchelor's Theorem~\cite{batchelor}, any super manifold~$N$ is isomorphic to $(|N|,\Gamma(\Lambda^\bullet E^*))$ for $E\to |N|$ a complex vector bundle over a smooth manifold~$|N|$. We denote such a super manifold by $\Pi E$. When doing geometry on super manifolds we will use the Deligne--Morgan sign conventions~\cite{strings1}; in particular, differential forms on super manifolds obey a bigraded sign rule. The super manifold~$\R^{n|m}$ is the locally ringed space with structure sheaf $C^\infty(\R^n)\otimes_\C\Lambda^\bullet(\C^m)$.

A \emph{vector bundle} over a super manifold is a finitely generated projective module over the structure sheaf. 
\emph{Generalized manifolds} and \emph{generalized super manifolds} are functors ${\sf Mfld}^{\rm op}\to {\sf Set}$ and ${\sf SMfld}^{\rm op}\to {\sf Set}$ respectively, i.e., presheaves on manifolds and super manifolds, respectively. We will use the notation $\underline{\sf Mfld}(M,N)$ to denote the generalized manifold $S\mapsto {\sf Mfld}(S\times M,N)$, and similarly for super manifolds, and we will frequently identify a super manifold with its representable functor. Generalized objects are \emph{representable} when they can be expressed as $S\mapsto {\sf Mfld}(S,M)$ for a manifold $M$. 
 For a (generalized) super manifold $\mathcal{M}$, we will use the notation $\mathcal{M}(S)$ to denote the set of maps $S\to \mathcal{M}$.  

\begin{ex}[The odd tangent bundle]\label{ex:oddT}
The generalized super manifold $\SM(\R^{0|1},M)$ plays a particularly important role in this paper. It turns out to be a representable super manifold, $\Pi TX\cong \SM(\R^{0|1},X)$, which identifies functions $C^\infty(\SM(\R^{0|1},X))\cong C^\infty(\Pi TX)\cong \Omega^\bullet(X)$ with differential forms on~$X$. The action of $\R^{0|1}$ on itself by translation gives an action map $\R^{0|1}\times \SM(\R^{0|1},X)\to \SM(\R^{0|1},X)$ whose derivative at $0\in \R^{0|1}$ is the de~Rham operator~$d$ viewed as an odd vector field on $\SM(\R^{0|1},X)$, i.e., an odd derivation on $C^\infty(\SM(\R^{0|1},X))\cong \Omega^\bullet(X)$. The $\Z/2$-action generated by reflection on~$\R^{0|1}$ also acts by precomposition, and is the pairity involution on $C^\infty(\SM(\R^{0|1},X))\cong \Omega^\bullet(X)$.
\end{ex}
 
 \subsubsection{Stacks}
For our purposes, smooth stacks are objects in the bicategory of (generalized) super Lie groupoids, bibundles and bibundle maps. Denote this bicategory by ${\sf SmSt}$. In particular, we identify super Lie groupoids with the stacks they present. An introduction to this perspective can be found in Section~2 of~\cite{BlohmannStacks}. In brief, a stack (up to isomorphism) is a Morita equivalence class of a Lie groupoid. 

Following the standard convention in geometry, our symmetry groups will always act on the left unless otherwise noted. An important consequence is that for a mapping space $M=\underline{\sf SMfld}(Y,X)$ a diffeomorphism $f\colon Y\to Y$ acts on $x\in M$ by $x\mapsto x\circ f^{-1}$; a diffeomorphism $g\colon X\to X$ will act by $x\mapsto g\circ x$. 

A \emph{quotient stack} is a stack that admits a presentation by a quotient groupoid, $M\sq G$, for $G$ acting on $M$. Equivariant vector bundles on the $G$-super manifold $M$ form a category equivalent to the category of vector bundles on the stack presented by $M\sq G$. In particular, a representation $\rho\colon G\to {\rm End}(V)$ defines a vector bundle $V_\rho$ on the stack presented by~$M\sq G$. When $V=\C$, sections of $V_\rho$ and $\Pi V_\rho$ are
\beq\begin{array}{ccc}
\Gamma(M\sq G, V_\rho)&\cong& \{f\in C^\infty(M) \mid \mu^*(f)=p_1^*(\rho)\cdot p_2^*(f) \in C^\infty(G\times M)\},\\
\Gamma(M\sq G,\Pi V_\rho)&\cong& \{f\in C^\infty(M)^{\odd}\mid \mu^*(f)=p_1^*(\rho)\cdot p_2^*(f)\in C^\infty(G\times M)\}\end{array}
\label{compsections}
\eeq
where $p_1\colon G\times M\to G,$ $p_2\colon G\times M\to M$ are the projection maps, and $\mu\colon G\times M\to M$ is the action map. Note that for the trivial representation sections are precisely the $G$-invariant functions on~$M$.

\begin{ex} The $\R^{0|1}\rtimes \Z/2$-action on $\SM(\R^{0|1},X)$ outlined in Example~\ref{ex:oddT} defines a quotient stack $\SM(\R^{0|1},X)\sq \R^{0|1}\rtimes \Z/2$. Functions on this stack are functions on $\SM(\R^{0|1},X)$ invariant under the action which are closed differential forms of even degree. The standard inclusion map $\Z/2=\{\pm 1\}\hookrightarrow \C^\times\cong {\rm GL}_1(\Pi \C)$ defines an odd line bundle over this stack whose sections are odd, closed differential forms.\label{ex:01stack}
\end{ex}

\subsubsection{Super determinant line bundles for families of invertible operators}\label{sec:sdet}

For a Fredholm operator~$D$ with discrete spectrum $\{\lambda_k\}_{k\in \Z}$, following Ray--Singer~\cite{RaySinger1,RaySinger2} we define the $\zeta$-function,
$$
\zeta_D(s)=\sum \lambda_k^{s}.
$$
In a wide range of examples, this defines a holomorphic function in $s$ for ${\rm Re}(s)\ll -1$ that can be analytically continued to a meromorphic function on~$\C$ that is regular at~$s=0$. We define the $\zeta$-determinant as
$$
{\rm det}_\zeta(D):=\exp(\zeta_D'(0)). 
$$
The $\zeta$-Pfaffian is a square root of the $\zeta$-determinant, and below we take~$\exp(\frac{1}{2}\zeta_D'(0))$. For operators that act on $\Z/2$-graded vector bundles, we define the \emph{$\zeta$-super determinant} as
$$
{\rm sdet}_\zeta(D):=\frac{{\rm pf}_\zeta(D|_{\rm odd})}{{\rm det}_\zeta(D|_{\rm even})^{1/2}}.
$$ 

The $\zeta$-determinant can also be applied to a family of operators parametrized by a smooth manifold~$M$, where each~$\lambda_k\in C^\infty(M)$ and ${\det}_\zeta(D)\in C^\infty(M)$. This procedure has an evident generalization to $M$ a super manifold by Taylor expanding~$\lambda_k^{-s}$ in odd variables, then computing the $s$-derivative at zero termwise in the Taylor series.

Using $\zeta$-determinant techniques, Quillen constructed a metrized determinant line bundle associated to twisted $\bar \partial$-operators on families of Riemann surfaces. This determinant line comes with a section~$\sigma$, and in cases where the bundle trivializes, $\sigma$ can be identified with the $\zeta$-determinant of~$D$ as a function. This setup was generalized by Bismut and Freed~\cite{BismutFreed1} for families of twisted Dirac operators. In brief, the fiber of the determinant bundle is the line $(\det {\rm ker}(D))^*\otimes (\det {\rm coker} (D))$ where~$\det$ denotes the highest exterior power. The subtlety is that the kernel and cokernel (although finite-dimensional) can jump in dimension as one moves around the family. However, for the purposes of this paper we only need to treat the case of \emph{invertible} families of operators, in which case the determinant line is canonically (topologically) trivialized. The section $\sigma$ in this case is the trivializing section, but the interesting feature is the metric is constructed so that
\beq
\|\sigma\|^2={\det}_\zeta(D^*D). \label{eq:detnorm}
\eeq
We take the obvious generalization for super families. 

\begin{defn} For a family of invertible operators~$D$ parametrized by a super manifold, the \emph{super determinant line} is the trivial line with trivializing section~$\sigma$ having norm squared~${\rm sdet}_\zeta(D^*D)$. \label{defn:supdet}
\end{defn}


Naturality of the determinant line construction under base change makes it amenable to stacky generalizations. This is done implicitly for the Dirac operator over the moduli stack of elliptic curves in \S4 of~\cite{Freed_Det}, wherein the Bismut--Freed holonomy formula~\cite{BismutFreed2} supplies the transformation properties for the determinant line under automorphisms of families of elliptic curves. Closely related transformation properties were also studied in detail by Atiyah~\cite{Atiyahlog}. 

In the case of a family of invertible operators over a (super) stack $\mathcal{M}$, the pullback of these operators to an $S$-point of~$\mathcal{M}$ leads to a determinant line bundle over $S$ that is canonically trivialized by its determinant section. However, this section can pick up nontrivial phases over isomorphisms between $S$-points. 
In general studying the phase of the determinant can be quite involved (e.g., $\eta$-invariant computations as in the Bismut--Freed holonomy formula), but in our cases there turns out to be an essentially unique way to fix the phases, gotten by descending trivializations at each $S$-point to the stack.

\subsubsection{Model super geometries}
Any super manifold $\M$ with a left action of a super group $G$ defines a \emph{model geometry}. Let an $(\M,G)$-super manifold consist of open submanifolds of $\M$ glued along (restrictions of) the action of $G$ on $\M$.\footnote{We require a cocycle condition when the $G$-action is not effective; see \cite{mingeo} Section 6.3 for details.}  Isometries of an $(\M,G)$-super manifold consist of diffeomorphisms that restrict locally to an action of $G$ on $\M$. For an $S$-family of super manifolds $F\to S$, denote by $\Iso_S(F)$ the group of \emph{isometries} over $S$ associated to the model geometry $(\M,G)$, i.e., isomorphisms $F\to F$ over $S$ that are isometries in the fibers. This defines a stack of $(\M,G)$-super manifolds on the site of super manifolds. See \cite{mingeo} Section~6.3 for details.

\begin{ex}[Euclidean spin geometries] \label{ex:EucModel} Let $\M=\R^d$ with its flat metric and spin structure, and take the isometry group $\Iso(\R^d)=\E^d\rtimes \Spin(d)$ where~$\E^d$ is the group of translations associated to~$\R^d$. Then $(\R^d,\E^d\rtimes \Spin(d))$ defines the \emph{Euclidean spin geometry}. 
\end{ex}

\begin{ex}[Super Euclidean geometries] \label{ex:superEucmod} Given data: (1)~a real vector space~$V$ with inner product; (2) a complex spinor representation~$\Delta$ of $\Spin(V)$; (3) a $\Spin(V)$-equivariant symmetric pairing $\Gamma\colon \Delta\otimes \Delta\to V_\C$ we define the super group
$$
(V\times \Pi \Delta)\times (V\times \Pi \Delta)\to (V\times \Pi \Delta),\quad (v,\sigma)\cdot (v',\sigma')=(v+v'+\Gamma(\sigma,\sigma'),\sigma+\sigma').
$$
When ${\rm dim}_\R(V)=d$ and ${\rm dim}_\C(\Delta)=\delta$, we write $\R^{d|\delta}:=V\times \Pi \Delta$ for the \emph{super Euclidean space} that carries an action of a group of \emph{super Euclidean translations} that we denote by~$\E^{d|\delta}$. There is an exact sequence of groups,
$$
0\to V \to (V\times \Pi \Delta)\to \Pi \Delta\to 0
$$
so that the ordinary translations of~$V$ form a subgroup of $\E^{d|\delta}$. We also have an action of $\Spin(V)$ on $V\times \Pi \Delta$ via the spinor representation on $\Delta$ and through the homomorphism $\Spin(V)\to \SO(V)$ on $V$. This defines a super group $\Iso(\R^{d|\delta}):=\E^{d|\delta}\rtimes \Spin(V)$, the \emph{super Euclidean isometry group.} The pair $\R^{d|\delta}$ and $\Iso(\R^{d|\delta})$ define a \emph{super Euclidean geometry}. \label{ex:modelsuper}
\end{ex}

\begin{notation}
To distinguish between translation groups and the super manifolds on which they act,~$\E^{d|\delta}$ will denote a group of super translations with underlying super manifold~$\R^{d|\delta}$. \label{not:supertrans}
\end{notation}

\begin{ex}[$0|1$-dimensional Euclidean geometry]\label{ex:01Euc} When $d=0$ and $\delta=1$, there is a unique pairing $\Gamma$, and we have $\Iso(\R^{0|1})\cong \E^{0|1}\rtimes \Z/2$. We can consider a groupoid whose objects are maps $\SM(\R^{0|1},X)$ from the Euclidean super point to~$X$ and whose morphisms are super Euclidean isometries between super points that are compatible with the map to~$X$. This groupoid has a presentation by the action groupoid, $\SM(\R^{0|1},X)\sq \Iso(\R^{0|1})$ and by the discussion in Examples~\ref{ex:oddT} and~\ref{ex:01stack}, we have an identification
$$
C^\infty\left(\SM(\R^{0|1},X)\sq \Iso(\R^{0|1})\right)\cong \Omega^{\ev}_{\rm cl}(X)
$$
between functions and closed even forms on~$X$. In~\cite{HKST}, the authors show that these functions define $0|1$-dimensional Euclidean field theories over~$X$. The present paper generalizes these ideas by replacing the groupoid~$\SM(\R^{0|1},X)\sq \Iso(\R^{0|1})$ with versions related to higher dimensional super Euclidean tori over~$X$. 
\end{ex}

\begin{rmk}\label{rmk:stack} The $0|1$-dimensional sigma model has as fields $\SM(\R^{0|1},X)$, which carries an action by the (isometry) group $\R^{0|1}\rtimes \Z/2$. The above example highlights the difference between working with the space versus a stack of fields: the former has as functions all differential forms, whereas the latter has as functions cocycles for de~Rham cohomology. 
\end{rmk}

\begin{ex} \label{ex:superfication}
Given a Euclidean spin manifold~$M$ in the sense of Example~\ref{ex:EucModel}, Stolz and Teichner define a \emph{superfication} functor (\cite{ST11}, Equation~4.14) with values in super Euclidean manifolds, given by
$$
\mathcal{S}(M):=\Pi (\Spin(M)\times_{\Spin(d)}\Delta),
$$
where $\Delta$ is the complex spinor representation defining the super Euclidean model geometry. That this is a functor comes from the homomorphism of isometry groups,
\beq
\E^d\rtimes \Spin(d)\hookrightarrow \E^{d|\delta}\rtimes \Spin(d)\cong \Iso(\E^{d|\delta})\label{eq:superficationgroup}
\eeq
where the source is the isometry group of $\R^d$ with its flat metric and spin structure. Hence super Euclidean manifolds are roughly spin manifolds with an enlarged isometry group. 
\end{ex}

\subsubsection{Weak modular forms}\label{sec:backMF}

A \emph{family of (2-dimensional) lattices} is a family of homomorphisms $\Lambda\colon S\times \Z^2\to \R^2\cong \C$ such that the ratio of the generators $\ell\colon S\times\{(1,0)\}\to \C$ and $\ell'\colon S\times\{(0,1)\}\to \C$ defines a map $\frac{\ell}{\ell'}\colon S\to \mathfrak{h}\subset \C$ with image in the upper half plane. Let $L$ denote the smooth manifold of these lattices; we have an evident diffeomorphism $L\cong \C^\times\times \mathfrak{h}$ that sends a pair of generators~$\ell,\ell'$ to $(\ell,\ell/\ell')\in (\C^\times\times \mathfrak{h})(S)$. There is an action of $\C^\times\times \SL_2(\Z)$ on~$L$ by
$$
\left(\mu,\left[\begin{array}{cc} a& b \\ c & d\end{array}\right],\ell,\ell'\right)\mapsto (\mu^2(a\ell+b\ell'),\mu^2(c\ell+d\ell')),\quad \mu\in \C^\times, \ \left[\begin{array}{cc} a& b \\ c & d\end{array}\right]\in \SL_2(\Z). 
$$

\begin{defn} \emph{Weak modular forms of weight $n/2$} are holomorphic functions $h$ on $L$ that are ${\rm SL}_2(\Z)$-invariant and have the property that $h(\mu\cdot \Lambda)=\mu^{-n}h(\Lambda)$ for $\mu\in \C^\times$. Taking products of holomorphic functions gives a graded ring, denoted ${\rm MF}$ whose degree $n$ piece, denoted ${\rm MF}_n$, are the weight $n/2$ weak modular forms. Define $\MF^n:=\MF_{-n}$. 
\end{defn}

The above is equivalent to the more common definition of modular forms of weight $n/2$ as holomorphic functions $h$ on the upper half plane~$\mathfrak{h}$ with the property $h\left(\frac{a\tau+b}{c\tau+d}\right)=\epsilon(\tau)^nh(\tau)$ for $\tau\in\mathfrak{h}$, $\left[\begin{array}{cc} a & b \\ c & d\end{array}\right]\in {\rm SL}_2(\Z)$ and $\epsilon(\tau)$ a holomorphic square root of $c\tau+d$. This comes from pulling back line bundles along the equivalence of stacks $\mathfrak{h}\sq {\rm MP}_2(\Z)\stackrel{\sim}{\to} L\sq \SL_2(\Z)\times \C^\times$ where ${\rm MP}_2(\Z)$ is the metaplectic group, the nontrivial $\Z/2$-extension of~$\SL_2(\Z)$. 

For $k>1$, the $2k^{\rm th}$ Eisenstein series is
$$
E_{2k}(\ell,\ell')=\sum_{n,m\in \Z_*^2} \frac{1}{(n\ell+m\ell')^{2k}}
$$
where $\Z^2_*$ denotes pairs $(m,n)\in \Z^2$, not both zero. The second Eisenstein series denoted~$E_2$ will refer to the holomorphic version, 
\beq\begin{array}{ccc}\displaystyle
E_2(\tau)&=&2\displaystyle\zeta(2)+\sum_{n\in \Z\setminus \{0\}}\sum_{m\in \Z} \frac{1}{(m+n\tau)^2},\\ E_2(\ell,\ell')&=&\displaystyle 2\zeta(2)/{\ell^2}+\sum_{n\in \Z\setminus \{0\}}\sum_{m\in \Z} \frac{1}{(m\ell+n\ell')^2},\end{array}\label{eq:2ndE}
\eeq
where we have given both the version as a function on $\mathfrak{h}$ and on $L$. 
It is not a modular form, but instead
$$
E_2\left(\frac{a\tau+b}{c\tau+d}\right)=(c\tau+d)^2E_2(\tau)-2\pi i c(c\tau+d),
$$ 
or in terms of a $\SL_2(\Z)$-invariant function on $L$ we have
$$
E_2(\mu\ell,\mu\ell')=\mu^{-2} E_2(\ell,\ell')+2\pi i\frac{(1-\bar\mu/\mu)}{\ell\bar\ell-\bar\ell\ell'}.
$$
Let $E^*_2$ denote the modular (but non-holomorphic) Eisenstein series, 
$$
E_2^*(\tau,\overline{\tau}):=\lim_{\epsilon\to 0^+} \left(\sum_{(m,n)\in \Z^2_*} \frac{1}{(m\tau+n)^2|m\tau+n|^\epsilon}\right).
$$
We have the relationship between the non-modular and non-holomorphic Eisenstein series,
\beq
E_2^*(\tau,\overline{\tau})=E_2(\tau)-\frac{\pi}{{\rm im}(\tau)},\quad\quad E_2^*(\ell,\overline{\ell},\ell'\bar\ell')=E_2(\ell,\ell')-\frac{2\pi i}{\ell\bar\ell'-\bar\ell\ell'}.\label{eq:2ndEisen}
\eeq
%

\subsection{Acknowledgements} It is a pleasure to thank Kevin Costello, Ryan Grady, Dmitri Pavlov and Stephan Stolz for helpful suggestions during the development of this work. I also thank my advisor, Peter Teichner, for his encouragement and insight, and Uli Bunke, Martin Olbermann, Stephan Stolz and an anonymous reviewer for their comments and suggestions that improved both the presentation and content. Lastly, I thank Owen Gwilliam for many invaluable discussions about perturbative quantum field theory. 

\section{Warm-up 1: classical vacua for the $1|1$-sigma model and ${\rm K}\otimes \C$}\label{sec:warmup1}

In this section we prove an analogous result to Theorem \ref{thethm} for ${\rm K}$-theory. We define a stack of classical vacua, denoted $\Phi^{1|1}_0(X)$, that consists of rigid conformal super circles equipped with a map to~$X$ that factors through an odd line. We construct a sequence of line bundles,~$\omega^{l/2}$, over $\Phi_0^{1|1}(X)$. As before, the assignment $X\mapsto \Gamma(\Phi_0^{1|1}(X);\omega^{\bullet/2})$ is a sheaf on the site of smooth manifolds. 

\begin{prop} There is a natural isomorphism of sheaves of graded graded algebras over~$\C$,
$$
\Gamma (\Phi_0^{1|1}(-);\omega^{\bullet/2})\stackrel{\sim}{\to} \left\{ \begin{array}{cc} \Omega^{\ev}_{\rm cl}(-) & \bullet={\rm even}\\ 
\Omega^{\odd}_{\rm cl}(-) & \bullet={\rm odd}\end{array}\right.
$$
whose target is the sheaf of 2-periodic closed differential forms. This realizes the source sheaf as a differential cocycle model for ${\rm K}\otimes \C$.
\label{prop:Kthy}
\end{prop}

The connection to K-theory comes from the de~Rham theorem and Chern isomorphism
$$
\Omega^{\ev/\odd}_{\rm cl}(X)\to {\rm HP}^{\ev/\odd}(X;\C)\cong {\rm K}^{\ev/\odd}(X)\otimes \C
$$
from 2-periodic closed forms to complexified K-theory. This implies that the isomorphism of sheaves in Theorem~\ref{prop:Kthy} identifies $\Gamma  (\Phi_0^{1|1}(X);\omega^{\bullet/2})$ as Hopkins--Singer differential cocycles for K-theory with complex coefficients. 
However, there are easier super geometric ways of cooking up differential forms as a $\Z/2$-graded algebra, e.g., Example~\ref{ex:01Euc}. The upshot of the approach taken below is that through its connection to super loop spaces, Proposition~\ref{prop:Kthy} permits several enhancements invisible from the perspective of de~Rham cohomology in parallel to the list after the statement of Theorem~\ref{thethm}. In fact, the progress made by F.~Han and F.~Dumutrescu on Stolz and Teichner's proposed model for ${\rm K}$-theory (over~$\Z$) as $1|1$-dimensional field theories allows us to say a bit more in this case. 
\begin{enumerate}
\item Constructions of Dumutrescu and Han define a functor from the category of vector bundles with connection on~$X$ to $C^\infty  (\Phi_0^{1|1}(X))$ (as a discrete category) via the \emph{super holonomy}. The result is the Chern form of the vector bundle, reviewed in~\S\ref{sec:chern}. 
\item For a finite group~$G$ acting on a manifold~$X$, sections of $\omega^{\bullet/2}$ over $\Phi^{1|1}_0(X\sq G)$ determine a non-Borel equivariant refinement of Proposition~\ref{prop:Kthy}, and (together with~(1) above) constructs the delocalized Chern character~\cite{DBEdiffK,DBE_EquivTMF}. This character simultaneously generalizes the ordinary Chern character (when $G=\{e\}$) and the character of a group representation~(when $X=\pt$). 
\item The isomorphism in Proposition~\ref{prop:Kthy} can be refined to a $\Z/2$-equivariant one for an involution on $\Phi_0^{1|1}(X)$ gotten from \emph{time-reversal} on the super circle. In \S\ref{sec:realunit} we identify concordance classes of $\Z/2$-invariant sections of $\omega^{\bullet/2}$ with~${{\rm KO}^\bullet(X)\otimes \C}$. 
\item A gerbe with connection on~$X$, denoted $\tau$, can be used to modify the line bundle to one denoted~$\omega^{\bullet/2+\tau}$. In \S\ref{subsec:twistedK} we identify concordance classes of sections of $\omega^{\bullet/2+\tau}$ with $\tau$-twisted ${\rm K}$-theory with complex coefficients. 
\item In Section~\ref{sec:warmup2} we construct a family of operators over $\Phi_0^{1|1}(X)$ whose $\zeta$-determinant is the $\hat{A}$-class of~$X$ as a function on~$\Phi_0^{1|1}(X)$. In this description, the appearance of the Bernoulli numbers (or equivalently, values of the Riemann $\zeta$-function) in the characteristic series for the $\hat{A}$-genus comes from their connection to traces of powers of the inverse of $id/dt$ acting on functions on $S^1$.
\end{enumerate}

In \S\ref{sec:super circles}-\ref{sec:prop:Kthy} we prove Proposition~\ref{prop:Kthy} before providing the physical motivation in~\S\ref{sec:SUSYmechinter} that guides the construction. We then explain the enhancements (1), (3), and (4) above in~\S\ref{sec:chern}, \S\ref{sec:realunit}, and \S\ref{subsec:twistedK}, respectively. 

\subsection{Super circles}\label{sec:super circles}


We begin by defining the $1|1$-dimensional rigid conformal model geometry as an extension of a super Euclidean model geometry of Example~\ref{ex:superEucmod}. 

\begin{defn} The \emph{$1|1$-dimensional super Euclidean model geometry} takes as data (1) $V=\R$ with its standard inner product, (2) the standard representation of~$\Spin(1)\cong \{\pm 1\}$ on~$\Delta=\C$, and (3) the nondegenerate pairing $\Gamma\colon \Delta\otimes\Delta \to V_\C$ given by $\Gamma(s,s')=iss'$. This determines a super group~$\E^{1|1}$ whose underlying super manifold is $\R^{1|1}$ with multiplication
$$
(t,\theta)\cdot (t',\theta')=(t+t'+i\theta\theta',\theta+\theta'), \quad (t,\theta),(t',\theta')\in \R^{1|1}(S).
$$
We obtain the semidirect product $\E^{1|1}\rtimes \Spin(1)$ from the involution $(t,\theta)\mapsto (t,\pm\theta)$, for $(t,\theta)\in\E^{1|1}(S)$ and $\Spin(1)\cong \{\pm 1\}$. Then the \emph{super Euclidean model space} is $\R^{1|1}$ with \emph{super Euclidean isometry group}~$\E^{1|1}\rtimes \Z/2$. The \emph{super rigid conformal} model geometry has the same model space, but isometry group $\Iso(\R^{1|1})=\E^{1|1}\rtimes \R^{\times}$ for the dilation action, $(t,\theta)\mapsto (\mu^2t,\mu\theta)$ by $\mu\in \R^\times(S)$ extending the action by $\Spin(1)\cong \{\pm 1\}\subset \R^\times$. 
\end{defn}

\begin{rmk} 
The standard super conformal model geometry has as isometries all diffeomorphisms of $\R^{1|1}$ that preserve the distribution generated by the vector field $\partial_\theta-i\theta\partial_t$. The rigid conformal isometry group is smaller. 
\end{rmk}


A \emph{family of (1-dimensional, oriented) lattices} is a family of homomorphisms $r\colon S\times \Z\to \R$ such that the image of the generators $S\times \{1\}$ is in the positive half-space, $\R_{>0}\subset \R$. Hence, the quotient for the $\Z$-action on $S\times \R$ is a family of metrized circles. We can promote this to a family of rigid conformal super circles using the inclusion $\R\subset \R^{1|1}$ of super manifolds and the homomorphism $\E\hookrightarrow \E^{1|1}$ where $\E=\R$ with the usual additive structure. Explicitly, for a family of lattices $r\colon S\times \Z\to \R$, we have an action map
$$
 S\times\R^{1|1}\times \Z\to S\times \R^{1|1},\quad (s,t,\theta,n)\mapsto (s,t+nr,\theta),\quad n\in \Z, (s,t,\theta)\in S\times\R^{1|1}.
$$
The quotient is a \emph{family of (rigid conformal) super circles} denoted~$S\times_r \R^{1|1}$. 

\begin{defn} The stack of \emph{rigid conformal super circles}, denoted $\mathcal{M}^{1|1}$, has objects over $S$ lattices $r\colon S\times \Z\to \R$ and morphisms rigid conformal isometries, $S\times_r \R^{1|1}\to S\times_{r'} \R^{1|1}.$ 
\end{defn}

\begin{lem} \label{lem:11supercirc} The stack $\mathcal{M}^{1|1}$ has a groupoid presentation with objects~$\R_{>0}$ and morphisms the bundle of groups~$(\E^{1|1}\rtimes \R^\times\times \R_{>0})/\Z$ that is a quotient by the $\Z$-action 
$$
\E^{1|1}\rtimes \R^\times\times \R_{>0}\times \Z\to \R^{1|1}\times \R^\times \times \R_{>0},\quad (t,\theta,\mu,r,n)\mapsto (t+nr,\theta,\mu).
$$
The source map is the projection, the target map is the action of $\R^\times$ on $\R_{>0}$ sending $(\mu,r)\mapsto \mu^2r$, and the unit map is induced by inclusion along the neutral element in~$\E^{1|1}\rtimes \R^\times$. 
\end{lem}


\begin{proof}[Proof of Lemma~\ref{lem:11supercirc}] The super manifold of objects comes from the map $S\cong S\times \{1\}\subset S\times \Z\to \R$ determining the family of lattices. It remains to understand the isometries. Locally they are the left action of $\E^{1|1}\rtimes \R^\times$ on $\R^{1|1}$, and lifting the action to the universal cover we observe that isometries are determined (possibly non-uniquely) by an $S$-point of $\E^{1|1}\rtimes \R^\times$. Since $\Z\subset \E\subset \E^{1|1}$ is in the center of the super translation group, an $S$-point of $\E^{1|1}\rtimes \R^\times$ induces the identity isometry when it is of the form $(nr,0,+1)\in (\E^{1|1}\rtimes \Z/2)(S)$. Finally, the action by dilations in $\R^\times$ on $\E^{1|1}$ preserves the image of $\E\subset \E^{1|1}$ and agrees with the square of the standard dilation action on this image, so changes the image of $\Z\subset \E\subset \E^{1|1}$ according to the claimed formula. This proves the lemma. 
\ep

\begin{defn} Define a complex line bundle $\omega^{1/2}$ over $\mathcal{M}^{1|1}$ determined by the map 
$$
(\E^{1|1}\rtimes \R^\times \times \R_{>0})/\Z\to \Z/2\cong \{\pm 1\}\subset \C^\times
$$
induced by projection and the sign homomorphism $\R^\times\cong \R_{>0}\times \{\pm 1\}\to \{\pm 1\}$, which determines a morphism of stacks $\mathcal{M}^{1|1}\to \pt\sq \C^\times$. 
\end{defn}



\subsection{Fields and classical vacua}

\begin{defn} The \emph{stack of fields}, denoted $\Phi^{1|1}(X)$, is the stack associated to the prestack whose objects over $S$ are pairs $(r,\phi)$ where $r\in \R_{>0}(S)$ determines a family of super circles $S\times_r \R^{1|1}$ and $\phi\colon S\times_r\R^{1|1}\to X$ is a map. Morphisms are commuting triangles
\beq
\begin{tikzpicture}[baseline=(basepoint)];
\node (A) at (0,0) {$S\times_r \R^{1|1}$};
\node (B) at (3,0) {$S\times_{r'} \R^{1|1}$};
\node (C) at (1.5,-1.5) {$X$};
\node (D) at (1.5,-.6) {$\circlearrowright$};
\draw[->] (A) to node [above=1pt] {$\cong$} (B);
\draw[->] (A) to node [left=1pt]{$\phi$} (C);
\draw[->] (B) to node [right=1pt]{$\phi'$} (C);
\path (0,-.75) coordinate (basepoint);
\end{tikzpicture}\label{11triangle}
\eeq
where the horizontal arrow is an isomorphism of $S$-families of rigid conformal $1|1$-manifolds. 

\end{defn}

\begin{defn}
The \emph{stack of classical vacua}, denoted $\Phi_0^{1|1}(X)$, is the full substack of $\Phi^{1|1}(X)$ generated by pairs $(r,\phi)$ where $\phi\colon S\times_r\R^{1|1}\to X$ factors through the map $\proj\colon S\times_r\R^{1|1}\to S\times\R^{0|1}$ induced by the projection $\R^{1|1}\to \R^{0|1}$. 
\end{defn}

\begin{rmk} One can view $\Phi_0^{1|1}(X)$ as super-analog of the constant loops in~$\Phi^{1|1}(X)$, or as those super loops invariant under the $S\times_r \E^1$-action on super circles. In~\S\ref{sec:SUSYmechinter} we explain how this is a stack of classical vacua in $1|1$-dimensional supersymmetric mechanics. 
\end{rmk}




\begin{prop} 
There is an equivalence of stacks
$$
\left( \begin{array}{c} (\E^{1|1}\rtimes \R^\times\times  \R_{>0})/\Z \times \SM(\R^{0|1},X)) \\ \downarrow \downarrow \\ \R_{>0}\times \SM(\R^{0|1},X)\end{array} \right) \stackrel{\sim}{\to} \Phi^{1|1}_0(X),
$$
where the quotient by $\Z$ on morphisms is the same is in Lemma~\ref{lem:11supercirc}, the source map is projection, and the target map is the projection to $\E^{0|1}\rtimes \R^\times \times \R_{>0}\times \SM(\R^{0|1},X)$ followed by the standard action of $\E^{0|1}\rtimes \R^\times$ on $\SM(\R^{0|1},X)$ and the $\R^\times$-action on $\R_{>0}$ as in Lemma~\ref{lem:11supercirc}. 
\label{prop:11present}
\end{prop}

\bp To verify the claim on objects, for a given $(r,\phi)$, we have $r\in \R_{>0}(S)$, and there is a unique $\phi_0\colon S\times \R^{0|1}$ such that $\phi=\phi_0\circ \proj$. Morphisms in $\Phi_0^{1|1}(X)$ are determined by isometries of circles which we computed in Lemma~\ref{lem:11supercirc}; this gives the super manifold of morphisms claimed in the proposition. 

It remains to compute the source and target maps. The source map is clearly the projection. For $(t,\theta,\mu)$ determining an isometry with source $(r,\phi)\in (\R_{>0}\times \SM(\R^{0|1},X))(S)$, the target is a family of super circles given by~$\mu^2r\in \R_{>0}(S)$, as in Lemma~\ref{lem:11supercirc}, and the isometry acts by precomposition on~$\phi$. Lifting to the universal cover, this is the action of $\E^{1|1}\rtimes \R^\times$ on~$\phi$ by precompositon, which acts on~$\phi_0$ through the homomorphism to $\E^{0|1}\rtimes \R^\times$, giving the claimed target map.
\ep


\subsection{The proof of Proposition \ref{prop:Kthy}}\label{sec:prop:Kthy}


\begin{proof}[Proof of Proposition~\ref{prop:Kthy}]  
By Proposition~\ref{prop:11present} and Equation~\ref{compsections}, sections of $\omega^{k/2}$ are functions on the objects, $\R_{>0}\times \SM(\R^{0|1},X)$, whose pullbacks along the source and target maps differ by~$\pm 1$, according to the parity of~$k$. Since the target map factors through an action of $\E^{0|1}\rtimes \R^\times\cong \E^{0|1}\rtimes (\Z/2\times \R_{>0})$ on $\R_{>0}\times \Pi TX$, we can identify sections as functions invariant under $\E^{0|1}\rtimes \R_{>0}$ and equivariant for the action of $\Z/2$. By the description of the action in Example~\ref{ex:oddT}, invariance under $\R^{0|1}$ requires the differential form component to be closed, and equivariance for the $\Z/2$-action fixes the parity of the form. The $\R_{>0}$-action is generated by the vector field 
$$
\beta:=2r\partial_r\otimes \id-\id\otimes \deg
$$ on $\R_{>0}\times \SM(\R^{0|1},X)$, where $\deg$ is the degree derivation on forms. Sections in the kernel of $\beta$ are spanned by those of the form $r^{k/2}\otimes \alpha$ for $\alpha\in\Omega^k_{\rm cl}(X)$. Hence, the map 
\beq
\Omega^{\bullet}_{\rm cl}(X)\stackrel{\sim}{\to} \Gamma (\Phi^{1|1}_0(X),\omega^{\bullet/2}),\quad \alpha\mapsto (2\pi r)^{k/2}\otimes \alpha,\ \alpha\in \Omega^k_{\rm cl}(X)\label{eq:11normpi}
\eeq
induces the claimed isomorphism of $\Z/2$ graded algebras. 

Naturality of~$\Phi_0^{1|1}(X)$ in~$X$ turns $\Gamma(\Phi_0^{1|1}(X);\omega^{\bullet/2})$ into a presheaf. We have shown that this presheaf is isomorphic to the sheaf of 2-periodic closed differential forms, so in fact $\Gamma(\Phi_0^{1|1}(X);\omega^{\bullet/2})$ is a sheaf, and the isomorphism is one of sheaves. 
\ep


The factor of $\pi$ in (\ref{eq:11normpi}) is a convenient normalization when comparing to the Chern form of a vector bundle, as we explain in \S\ref{sec:chern}.

\subsection{Motivation from $1|1$-dimensional supersymmetric mechanics}\label{sec:SUSYmechinter}

Physical intuition for $\Phi_0^{1|1}(X)$ comes from the $1|1$-dimensional sigma model with target~$X$, which is supersymmetric mechanics  with minimal supersymmetry. A mathematical introduction can be found in D. Freed's notes~\cite{5lectures}, Chapter~2 and some computations are worked out in detail in~\cite{strings1} pages 649-656. We work with a Wick-rotated version of this theory; our conventions for this rotation, its affect on the Euclidean group, and its affect on the classical action functional are taken from P.~Deligne and D.~Freed's article, ``Classical Field Theory" in~\cite{strings1}.

\begin{defn} Let $S\times_r\R$ denote a family of loops equipped with a fiberwise spin structure. An \emph{$S$-family of fields for the $1|1$-dimensional sigma model with target~$X$} consists of: (1) a map~$x\colon S\times_r \R \to X$ and (2) an odd spinor $\psi\in \Gamma(S\times_r \R,\Pi \$\otimes x^*TX)$ twisted by the pullback tangent bundle. 
For a choice of metric on~$X$, define the \emph{classical action} 
\beq\begin{array}{ccl}
\mathcal{S}(x,\psi)
&=&\displaystyle\frac{i}{2}\int_0^r \left(\langle \dot{x},\dot{x}\rangle-i\langle \psi,\nabla_{\dot{x}}\psi\rangle\right)dt
\end{array}\label{eq:11classactcomp}
\eeq
where the integral is understood as being fiberwise over $S$. 
\end{defn}

The action $\mathcal{S}$ defines a function on the generalized super manifold whose $S$-points are $S$-families of fields. The Euler--Lagrange equations for this classical action are (\cite{strings1} page~654, or~\S\ref{sec:11Hess})
$$
\nabla_{\dot x} \dot x=-\frac{i}{2} R(\psi,\psi)\dot x,\quad\quad \nabla_{\dot x}\psi=0.
$$
Fields satisfying these equations are \emph{classical solutions}, and form a generalized sub-super manifold of the fields. The \emph{energy} of a classical solution is the fiberwise integral of $\frac{1}{2}|\dot x|^2$ on~$S\times_r \R$.

\begin{defn} A \emph{classical vacuum} over $S$ is an energy zero family of classical solutions. \end{defn}

\begin{lem} An $S$-family of fields $(x,\psi)$ is a classical vacuum if and only if~$x$ is a family of constant maps and~$\psi$ is a constant section. \end{lem}

\bp If the energy of a field is zero, then $\dot x=0$, and so $x\colon S\times_r \R\to X$ factors through the projection $S\times_r \R\to S$ as claimed. By virtue of being a classical solution, $\psi$ is covariantly constant, i.e., a constant section. \ep

We can repackage the data of the fields as maps from a family of super circles,~$S\times_r\R^{1|1}$. The lift of the $r\Z$-action on $S\times \R$ to $S\times \R^{1|1}$ depends on the fiberwise spin structure. In the odd (or nonbounding) case, $S$-families of fields are precisely objects of $\Phi^{1|1}(X)$ over~$S$, as the odd spinor bundle is a trivial bundle over $S\times_r \R$. We verify that a $S$-family of vacua in this case is a map $S\times_r \R^{1|1}\to X$ factoring through $S\times\R^{0|1}$: a trivialization of the spinor bundle identifies the covariantly constant section~$\psi$ with a map $S\times \R^{0|1}\to X$. 

So far in this subsection, we have dealt with fields as a sheaf of sets on the site of super manifolds, with the action being a function i.e., natural transformation to~$\C$. However, a classical field theory has a bit more structure: symmetries allow us to view the classical action as a function on a \emph{stack}. These symmetries are precisely the super Euclidean isometries between families of super circles. Extending from the super Euclidean isometries to rigid conformal ones includes the action of dilations of circles (via $\R^\times\cong \R_{>0}\times \Z/2$), which is precisely the action of the renormalization group (RG) flow. In summary, the stack~$\Phi_0^{1|1}(X)$ has as objects over~$S$ classical vacua for super circles with odd spin structures, and morphisms over~$S$ are symmetries of the classical action and the RG-flow. 

Although the classical action depended on a choice of metric on~$X$, the stack $\Phi_0^{1|1}(X)$ is independent of this choice. This is a feature that is special to the supersymmetric sigma models considered in this paper---in general, classical vacua will depend on the input geometric structures. 

\subsection{Differential cocycles from vector bundles} \label{sec:chern} 

Let $1|1\EFT(X)$ denote the category of $1|1$-dimensional (untwisted) field theories, defined similarly to the discussion in \S\ref{sec:Stolz--Teichner}. Let ${\sf Vect}^\nabla(X)$ denote the groupoid of vector bundles with connection on $X$ and connection preserving bundle automorphisms. We will explain the commutative diagram
\beq
\begin{tikzpicture}[baseline=(basepoint)];
\node (A) at (0,0) {${\sf Vect}^\nabla(X)$};
\node (B) at (3,0) {$1|1\EFT(X)$};
\node (C) at (0,-1.5) {$K^0(X)$};
\node (D) at (3,-1.5) {${\rm H}_{\rm dR}^{\ev}(X)$}; 
\node (E) at (7,-.75) {$\Gamma (\Phi_0^{1|1}(X);\omega^{\bullet/2})$};
\draw[->] (A) to node [above=1pt] {${\rm sPar}$} (B);
\draw[->] (A) to (C);
\draw[->] (C) to node [above=1pt] {${\rm ch}$} (D);
\draw[->] (B) to node [right=1pt] {${\rm ch}_{\rm Han}$} (D);
\draw[->] (B) to node [above=1pt] {$\partition$} (E);
\draw[->] (E) to node [below=1pt] {${\rm conc}$} (D);
\path (0,-.75) coordinate (basepoint);
\end{tikzpicture}\nonumber
\eeq
where ${\rm ch}$ is the Chern character, ${\rm ch}_{\rm Han}$ is the \emph{Han character} constructed in his thesis~\cite{Han}, $\partition$ evaluates a $1|1$-Euclidean field theory on $\Phi^{1|1}_0(X)$ (analogous to the functor $Z$ in~\S\ref{sec:Stolz--Teichner}), ${\rm conc}$ applies the inverse to the isomorphism (\ref{eq:11normpi}) and takes concordance classes, and ${\rm sPar}$ is the functor constructed by F.~Dumitrescu in his thesis~\cite{florin_11} that produces a super Euclidean field theory from a vector bundle with connection using super parallel transport. We remark that it remains an open question as to whether ${\rm sPar}$ surjects onto concordance classes of field theories. 

For a $1|1$-EFT in the image of the functor~${\rm sPar}$, evaluation on $\Phi_0^{1|1}(X)$ yields the function~${\rm Tr}(\exp(ir F))$ on $R_{>0}\times \Pi TX$ for $F$ the curvature 2-form of the vector bundle and $r$ the usual coordinate on~$\R_{>0}$. This was first computed by Han, though he actually produces the Bismut--Chern character, a richer but more complicated construction~\cite{Han}. We'll sketch an abbreviated version for the usual Chern character. 

For $E\to X$ a vector bundle with connection and $\phi\colon S\times \R^{1|1}\to X$ a super path factoring through $S\times \R^{0|1}$, define an operator~$\nabla_D$ on sections $s\in \Gamma(S\times \R^{1|1},\phi^*E)$ from the pullback connection, where $D=\partial_\theta-i\theta\partial_t$. Take as component fields $s_0=i_0^*s$ and $s_1=i_0^*\nabla_Ds$ for $i_0\colon S\times \R\hookrightarrow S\times \R^{1|1}$ the inclusion of the reduced manifold. Then the super parallel transport equation, $\nabla_Ds=0$, in terms of the component fields reads
$$
\nabla_{i\partial_t} s_0=\frac{1}{2} i_0^*(\phi^*F(D,D))s_0,\quad\quad s_1= 0,
$$
where we have used 
$$
\nabla_D\nabla_D=\frac{1}{2}(\nabla_D\nabla_D+\nabla_D\nabla_D)=\frac{1}{2}(F(D,D)-\nabla_{[D,D]})=\frac{1}{2}F(D,D)+\nabla_{i\partial_t}
$$ 
Hence, the solution to these differential equations for an initial condition $s_0$ is
$$
s_0(t)=\exp(it\cdot  i_0^*(\phi^*F(D,D))/2)s_0(0),\quad \quad s_1(t)\equiv 0.
$$
A quick computation (e.g., in local coordinates) finds that trace at $t=r$ gives the claimed formula in the universal case $S=\R_{>0}\times \Pi TX$, and~$\phi$ is 
$$
(\R^{1|1}\times \R_{>0})/\Z\times \Pi TX\to \R^{0|1}\times\Pi TX\cong \R^{0|1}\times \SM(\R^{0|1},X)\stackrel{\ev}{\to} X
$$ 
where the first arrow is the projection and $\ev$ is evaluation. 

\subsection{Time reversal, conjugation and ${\rm KO}$-theory}\label{sec:realunit}

Our model for ${\rm K}\otimes \C$ has a simple and geometrically-motivated homotopy $\Z/2$-action that is the complexification of the usual $\Z/2$-action on complex K-theory with homotopy fixed points ${\rm KO}$. 

\emph{Time reversal} is an automorphism of the model space $\R^{1|1}$ and isometry group $\E^{1|1}\rtimes \R^\times$ determined by~$(t,\theta)\mapsto (-t,i\theta)$. This acts on the charts of a $1|1$-rigid conformal manifold by precomposition and on gluing maps by conjugation, so determines an automorphism of the stack of $1|1$-rigid conformal manifolds. This determines an endofunctor ${\sf t}\colon \Phi_0^{1|1}(X)\to \Phi_0^{1|1}(X)$ with ${\sf t}\circ {\sf t}\cong \id$, i.e., a homotopy $\Z/2$-action. 

On $C^\infty(\R_{>0}\times \Pi TX)$, the action by ${\sf t}$ is $i^{\deg}$ on differential forms. In particular, functions that are invariant under time-reversal are the 4-periodic forms, so we obtain cocycles,
$$
C^\infty(\Phi_0^{1|1}(X))^{\sf t}\to \bigoplus_{i} \Omega^{4i}_{\rm cl}(X)\twoheadrightarrow {\rm KO}(X)\otimes \C. 
$$ 
In fact, over~$\C$ the action on coefficients determines the action completely, so the ring $C^\infty (\Phi_0^{1|1}(X))$ with its ${\sf t}$-action is a cocycle model for ${\rm KR}(X)\otimes \C$ with its $\Z/2$-action.

\subsection{Gerbes and twisted ${\rm K}$-theory}\label{subsec:twistedK}

Bundle gerbes with connection provide a geometric model for twists of ${\rm K}$-theory. For a gerbe with curvature 3-form $H$, twisted ${\rm K}$-theory with complex coefficients can be computed by the de~Rham complex with the twisted differential~$d+H$. To do this in a manner which can be made functorial, for a cover $\{U_i\}$ of $X$, a class in twisted cohomology consists of forms $\alpha'_i$ on $U_i$ such that $\alpha'_j=\exp(F_{ij})\alpha'_i$ on $U_i\cap U_j$ where $F_{ij}$ denotes the curvature of the line bundle with connection defined on $U_i\cap U_j$ in the cocycle data for a gerbe with connection on $X$; e.g., see Atiyah and Segal's account in \cite{AtiyahSegaltwists}, particularly Section 6. In this subsection we show that a gerbe with connection on $X$ determines a line bundle over classical vacua whose sections are classes in twisted ${\rm K}$-theory with complex coefficients.

Let cocycle data for a bundle gerbe with connection on $X$ be given, i.e., suppose we have a cover $\{U_i\}$ of $X$ together with line bundles with connection $( L_{ij},\nabla_{ij})$ for each overlap $U_i\cap U_j$ satisfying a cocycle condition. These line bundles determine super parallel transport functors for super paths in $U_i\cap U_j$, and restriction to $\Phi^{1|1}_0(X)$ gives a function on $\R_{>0}\times \Pi T (U_i\cap U_j)$ coming from the holonomy around these super paths; when restricted to $\R_{>0}\times \Pi T (U_i\cap U_j)$ this is precisely the function $\exp(irF_{ij})$ for $F_{ij}$ the curvature of the connection $\nabla_{ij}$. These nonvanishing functions on overlaps determine a line bundle on $\R_{>0} \times \Pi TX$, where the cocycle condition can be verified immediately via the assumed cocycle condition for the gerbe with connection. 

The vector space $\Gamma (\Phi_0^{1|1}(X);\omega^{l/2+\tau})$ is spanned by functions $\alpha'_i\in \Gamma (\Phi_0^{1|1}(U_i);\omega^{l/2})$ on each $U_i$ such that $\alpha'_j=\exp(irF)\alpha'_i$. Hence, we've shown that a gerbe with connection leads to a natural isomorphism of algebras,
$$
{\rm K}^{\bullet+\tau}(X)\otimes \C\cong \Gamma (\Phi_0^{1|1}(X);\omega^{\bullet/2+\tau})/{\rm concordance}
$$
where the left hand side denotes ${\rm K}$-theory twisted by the gerbe, taken with complex coefficients, and the right hand side denotes sections of the line bundle over classical vacua constructed out of $\tau$. 

\section{Classical vacua of the $2|1$-sigma model and ${\rm TMF}\otimes \C$}

Now we pass to the $2|1$-dimensional case, defining a stack of classical vacua denoted $\Phi^{2|1}_0(X)$ and a line bundle $\omega^{1/2}$ over this stack. With these definitions in place, proving Theorem~\ref{thethm} is a similar computation to the $1|1$-dimensional version. 

\subsection{Super tori}

We begin with the $2|1$-dimensional rigid conformal model geometry (called \emph{conformal Euclidean} in~\cite{ST11}) which is an extension of the super Euclidean model geometry. 

\begin{defn} The \emph{$2|1$-dimensional Euclidean model geometry} takes as data (1) $V=\R^2\cong \C$ with its standard inner product, (2) the dual of the standard representation of~$\Spin(2)\cong U(1)$ on~$\Delta=\C$, and (3) the nondegenerate pairing $\Gamma\colon \Delta\otimes\Delta \to V_\C\cong \C_1\oplus \C_{-1}$ that maps to $\C_{-1}$, where the subscript labels the action of $\SO(2)\cong U(1)$ on the subspace of~$V_\C$. This determines a super group~$\E^{2|1}$ of translations whose underlying super manifold is $\R^{2|1}$ with multiplication
$$
(z,\bar{z},\theta)\cdot (z',\bar{z}',\theta')=(z+z',\bar{z}+\bar{z}'+\theta\theta',\theta+\theta'),\quad (z,\bar{z},\theta),\ (z',\bar{z}',\theta')\in \R^{2|1}(S).
$$
The \emph{super Euclidean group} is the semidirect product $\E^{2|1}\rtimes \Spin(2)$. The \emph{rigid conformal} model geometry has the same model space with isometry group 
$$
\Iso(\R^{2|1})=\E^{2|1}\rtimes \C^\times\cong \E^{2|1}\rtimes (U(1)\times \R_{>0})\cong \E^{2|1}\rtimes (\Spin(2)\times \R_{>0})
$$
with the action of $\C^\times$ by 
$$
(\mu,\bar{\mu})\cdot(z,\bar{z},\theta)=(\mu^2z,\bar{\mu}^2\bar{z},\bar{\mu}\theta),\quad (\mu,\bar{\mu})\in \C^\times(S),\ (z,\bar{z},\theta)\in \R^{2|1}(S).
$$
\end{defn}

In \S\ref{sec:backMF}, we defined a family of (based, oriented) lattices as a family of homomorphisms $\Lambda\colon S\times \Z^2\to \R^2\cong \C$ such that the ratio of the generators $\ell\colon S\times\{(0,1)\}\to \C$ and $\ell'\colon S\times\{(1,0)\}\to \C$ defines a map $\frac{\ell}{\ell'}\colon S\to \mathfrak{h}\subset \C$ with image in the upper half plane. In particular, a family of lattices defines an $S$-family of tori we denote by~$S\times_\Lambda \R^2$ gotten by the quotient of $S\times \R^2$ by the fiberwise $\Z^2$-action. The inclusion $\R^2\subset \R^{2|1}$ of super manifolds and $\C\cong \E^2\subset \E^{2|1}$ of super groups allows us to define a family of super tori for any lattice via the quotient $S\times \R^{2|1}$ by the action of $\Lambda$. We denote this quotient by~$S\times_\Lambda \R^{2|1}$. 

\begin{defn} The stack of \emph{(super) rigid conformal tori}, denoted $\mathcal{M}^{2|1}$, has objects over~$S$ lattices $\Lambda\colon S\times \Z^2\to \C$ and morphisms rigid conformal isometries, $S\times_\Lambda \R^{2|1}\to S\times_{\Lambda'} \R^{2|1}.$ 
\end{defn}

\begin{lem} \label{lem:21supercirc} The stack $\mathcal{M}^{2|1}$ has a groupoid presentation with objects~$L$ and morphisms the quotient $(\E^{2|1}\rtimes \C^\times \times \SL_2(\Z) \times L)/\Z^2$ for the $\Z^2$-action 
\beq
\E^{2|1}\rtimes \C^\times \times\SL_2(\Z)\times L\times \Z^2&\to& \E^{2|1}\rtimes \C^\times\times \SL_2(\Z) \times L,\nonumber\\
 (z,\bar z,\theta,\mu,A,\ell,\bar \ell,\ell'\bar\ell',n,n')&\mapsto& (z+n\ell+n'\ell',\bar z+n\bar\ell+n'\bar\ell',\theta,\mu,A,\ell,\bar\ell,\ell',\bar\ell').\nonumber
\eeq
for
$$
(n,n')\in \Z^2(S), \ (z,\bar z,\theta)\in \E^{2|1}(S), \ \mu\in \C^\times(S),\ A\in \SL_2(\Z)(S), \ (\ell,\bar \ell,\ell'\bar\ell')\in L(S).
$$
The source map is the projection, the target map is the $\C^\times\times \SL_2(\Z)$-action on~$L$, and the unit map is induced by inclusion along the neutral element in~$\E^{2|1}\rtimes \C^\times\times \SL_2(\Z)$. 
\end{lem}


\begin{proof}[Proof of Lemma~\ref{lem:21supercirc}] The objects are $L$ by definition, so it remains to understand isometries between $S$-family of rigid conformal tori,
$$
S\times_\Lambda\R^{2|1}\to S\times_{\Lambda'}\R^{2|1}.
$$ 
One difference from the $1|1$-dimensional case is that the lattices $\Lambda$ and $\Lambda'$ can both differ from a dilation \emph{and} an $S$-point of $\SL_2(\Z)$. Since the (left) super translation action commutes with the (right) lattice action on $\R^{2|1}$, this verifies the source and target maps are as claimed, up to understanding the group structure on the super translations themselves. These remaining isomorphisms are locally the left action of $\E^{2|1}$ on $\R^{2|1}$, and lifting the action to the universal cover we observe that action by super translations is determined (possibly non-uniquely) by an $S$-point of $\E^{2|1}$. Since $\Z^2\subset \E\subset \E^{2|1}$ is in the center of the super translation group, an $S$-point of $\E^{2|1}$ induces the identity isometry when it is a translation in the lattice $S\times \Z^2\stackrel{\Lambda}{\hookrightarrow} S\times \E^{2|1}.$ Hence, the morphism space is the claimed quotient. 
\ep

\begin{defn} Define a complex line bundle $\omega^{1/2}$ over $\mathcal{M}^{1|1}$ determined by the projection 
$$
(\E^{2|1}\rtimes \C^\times\times \SL_2(\Z)\times L)/\Z^2\to \C^\times
$$
that induces a morphism of stacks $\mathcal{M}^{2|1}\to \pt\sq \C^\times$. \end{defn}

Many quantities below depend on the volume of a torus,
\beq
\vol:=\frac{1}{2i}(\ell\bar\ell'-\bar\ell\ell')\in C^\infty(L).\label{eq:vol}
\eeq
Note that when $(\ell,\ell')=(\tau,1)$, we get $\vol={\rm im}(\tau)$, as usual. 


\subsection{Fields and classical vacua}

\begin{defn} The \emph{stack of fields}, denoted $\Phi^{2|1}(X)$, is the stack associated to the prestack whose objects over $S$ are pairs $(\Lambda,\phi)$ where $\Lambda\in  L(S)$ determines a family of super tori $S\times_\Lambda \R^{2|1}$ and $\phi\colon S\times_\Lambda \R^{2|1}\to X$ is a map. Morphisms over~$S$ are commuting triangles
\beq
\begin{tikzpicture}[baseline=(basepoint)];
\node (A) at (0,0) {$S\times_\Lambda \R^{2|1}$};
\node (B) at (3,0) {$S\times_{\Lambda'} \R^{2|1}$};
\node (C) at (1.5,-1.5) {$X$};
\node (D) at (1.5,-.6) {$\circlearrowright$};
\draw[->] (A) to node [above=1pt] {$\cong$} (B);
\draw[->] (A) to node [left=1pt]{$\phi$} (C);
\draw[->] (B) to node [right=1pt]{$\phi'$} (C);
\path (0,-.75) coordinate (basepoint);
\end{tikzpicture}\label{21triangle}
\eeq
where the horizontal arrow is an isomorphism of $S$-families of super rigid conformal $2|1$-dimensional super manifolds. 

\end{defn}

\begin{defn} 
The \emph{stack of classical vacua}, denoted $\Phi_0^{2|1}(X)$, is the full substack of $\Phi^{2|1}(X)$ generated by pairs $(\Lambda,\phi)$ where $\phi\colon S\times_\Lambda\R^{2|1}\to X$ is a map that factors through the map $\proj\colon S\times_\Lambda\R^{2|1}\to S\times \R^{0|1}$ induced by the projection $\R^{2|1}\to \R^{0|1}$. 
\end{defn}

\begin{rmk} The inclusion $\Phi_0^{2|1}(X)\subset \Phi^{2|1}(X)$ is a super-geometric analog of the constant maps of tori including into the free double loop space, or equivalently, the inclusion of maps invariant under the action of the reduced tori~$S\times_\Lambda \E^2$. In~\S\ref{sec:21motivate} we explain how $\Phi_0^{2|1}(X)$ can be viewed as the classical vacua for the $2|1$-dimensional sigma model. \end{rmk}

\begin{prop} 
There is an equivalence of stacks
$$
\left( \begin{array}{c} (\E^{2|1}\rtimes \C^\times\times\SL_2(\Z)\times L)/\Z^2 \times \SM(\R^{0|1},X)) \\ \downarrow \downarrow \\ L\times \SM(\R^{0|1},X)\end{array} \right) \stackrel{\sim}{\to} \Phi^{2|1}_0(X),
$$
where the quotient by $\Z^2$ on morphisms is the same is in Lemma~\ref{lem:21supercirc}, the source map is projection, and the target map is the projection to $\E^{0|1}\rtimes \C^\times\times \SL_2(\Z)\times L\times \SM(\R^{0|1},X)$ followed by the action of $\E^{0|1}\rtimes \C^\times$-action on $\SM(\R^{0|1},X)$ and the $\C^\times \times \SL_2(\Z)$-action on~$L$.
\label{prop:21present}
\end{prop}

\bp To verify the claim on objects, for a given $(\Lambda,\phi)$, we have $\Lambda \in L(S)$, and there is a unique $\phi_0\colon S\times \R^{0|1}$ such that $\phi=\phi_0\circ \proj$. Morphisms in $\Phi_0^{2|1}(X)$ are determined by isometries of circles which we computed in Lemma~\ref{lem:21supercirc}; this gives the super manifold of morphisms claimed in the proposition. 

It remains to compute the source and target maps. The identity isometry leaves~$(\Lambda,\phi)$ unchanged, so the source map is the projection. Isometries from $\SL_2(\Z)$ act on~$L$, and this acts trivially on the map $\phi$; those from $\E^{2|1}\rtimes \C^\times$ act on $\Lambda$ through the projection to $\C^\times$, and on the map $\phi$ through precomposition with an isometry of a family of super circles. Lifting to the universal cover, this is determined by an action of~$\E^{2|1}\rtimes \C^\times$ on~$\phi$, which acts on~$\phi_0$ through the homomorphism to $\E^{0|1}\rtimes \C^\times$. This verifies the claimed target map.
\ep

\subsection{The proof of Theorem~\ref{thethm}}\label{sec:proofthm1}

We have $\mathcal{M}^{2|1}\cong \Phi^{2|1}(\pt)\cong \Phi_0^{2|1}(\pt)$, so $(\omega^{1/2})^{\otimes k}=:\omega^{\bullet/2}$ also defines line bundles on these stacks. The canonical map $X\to \pt$ induces a functor $\Phi^{2|1}_0(X)\to \Phi^{2|1}_0(\pt)$ along which we can pull back $\omega^{\bullet/2}$. First we consider the space of all sections before defining and characterizing the holomorphic ones. 



\begin{prop} There is an isomorphism of graded algebras, 
$$
\bigoplus_{i+j=k} \Omega^j_{\rm cl}(X;\MaF^i)\stackrel{\sim}{\hookrightarrow} \Gamma(\Phi_0^{2|1}(X),\omega^{k/2}),\quad \alpha\otimes F\mapsto \vol^{j/2}\alpha\otimes F,
$$
where $\alpha\in \Omega^j_{\rm cl}(X)$ and $F\in \MaF^j$, and we regard $\vol^{j/2}\alpha\otimes F$ as a function on the cover $L\times \Pi TX\cong L\times \SM(\R^{0|1},X)$ that descends to a section over the stack $\Phi_0^{2|1}(X)$. 
  \label{prop:Maass}
\end{prop}

\begin{proof}[Proof of Proposition~\ref{prop:Kthy}] By Proposition~\ref{prop:21present} and Equation~\ref{compsections}, sections of $\omega^{k/2}$ are functions on the objects, $L\times \SM(\R^{0|1},X)\cong L\times \Pi TX$ with equivariant properties when pulled back to $S$-families of isomorphisms between objects. It will be convenient to express a function $f\in C^\infty(L)\otimes \Omega^j(X)\subset C^\infty(L\times \Pi TX)$ as
$$
f=\vol^{j/2} F\otimes \alpha,\quad F\in C^\infty(L),\ \alpha\in \Omega^j(X),
$$
for $\vol$ the volume of a torus defined in (\ref{eq:vol}). Now, equivariance under isometries implies in particular that $f$ is invariant under $\SL_2(\Z)\times \E^{2|1}$. 
The volume of a torus is an $\SL_2(\Z)$-invariant function, so~$f$ is $\SL_2(\Z)$-invariant if and only if~$F\in C^\infty(L)^{\SL_2(\Z)}$. The $\R^{2|1}$-action is through the de~Rham $d$ and so $\alpha\in \Omega^j_{\rm cl}(X)$ is closed. Hence, $f\in C^\infty(L)^{\SL_2(\Z)}\otimes \Omega^j_{\rm cl}(X)$. 

We claim that equivariance of $f$ for the $\C^\times$-action makes~$F$ transform with weight $(k-j)/2$. The dilation action on a $j$-form is through~$\bar{\mu}^{-j}$, so 
\beq
(\mu,\bar{\mu})\cdot \left(\vol^{j/2}F\otimes \alpha \right)&=&(\mu^2\bar\mu^2 \vol)^{j/2} \left((\mu,\bar \mu)\cdot F\right)\otimes \bar{\mu}^{-j}\alpha\nonumber\\
&=&\mu^j\vol^{j/2}\left((\mu,\bar{\mu})\cdot F\right)\otimes \alpha,\nonumber
\eeq
so for $\vol^{j/2} F\otimes \alpha$ to be a section of $\omega^{k/2}$, we require that $(\mu,\bar{\mu})\cdot F=\mu^{k-j} F$. Hence, sections of $\omega^{k/2}$ are generated by sums of functions $\vol^{j/2}F \otimes \alpha$ for $F\in {\rm MaF}^i$ and $\alpha\in \Omega^j_{\rm cl}(X)$ such that $i+j=k$. We choose the isomorphism that takes a homogeneous element $F\otimes \alpha \in \Omega^j_{\rm cl}(X;\MaF^i)$ to $\vol^{j/2}F \otimes \alpha\in \Gamma(\Phi_0^{2|1}(X);\omega^{k/2})$. 

Naturality of~$\Phi_0^{2|1}(X)$ in~$X$ turns $\Gamma(\Phi_0^{2|1}(X);\omega^{\bullet/2})$ into a presheaf. We have shown that this presheaf is isomorphic to the sheaf of closed differential forms with values in weak Maass forms, so in fact $\Gamma(\Phi_0^{2|1}(X);\omega^{\bullet/2})$ is a sheaf, and the isomorphism is one of sheaves. 
\ep

Using the same notation as in the proof above, the dilation equivariance property for~$F$ under the $\C^\times$-action has an infinitesimal expression in terms of the generators of the Lie algebra of $\C^\times$, namely
$$
2(\ell\partial_\ell+\ell'\partial_{\ell'})F=(k-j)F,\quad\quad 2(\bar\ell\partial_{\bar\ell}+\bar\ell'\partial_{\bar\ell'})F=0,
$$
where, e.g., the vector field $\bar\ell\partial_{\bar \ell}$ corresponds to dilation in the variable~$\bar\ell$ on~$L$, or more precisely, the anti-holomorphic part of the vector field that dilates~$(\ell,\bar\ell)$. We can demand invariance under each of the infinitesimal anti-holomorphic dilations separately, 
$$
2\bar\ell\partial_{\bar\ell}F=0,\quad\quad\quad 2\bar\ell'\partial_{\bar\ell'} F=0
$$
so that $F$ is holomorphic (in the ordinary sense) as a function on $L\cong \mathfrak{h}\times \C^\times$. Tracing through the computation in the previous proof, these conditions on~$F$ are equivalent to
\beq
\left(2\bar\ell \partial_{\bar \ell}-\frac{\bar\ell\ell'}{\bar\ell\ell'-\bar\ell'\ell} \deg\right)f=0,\quad \quad \left(2\bar\ell' \partial_{\bar \ell'}+\frac{\bar\ell'\ell}{\bar\ell\ell'-\bar\ell'\ell} \deg\right)f=0\label{eq:holo}
\eeq
for $f\in \Gamma(\Phi_0^{2|1}(X);\omega^{k/2})$, where $\deg$ is the degree derivation that coincides with the infinitesimal action on functions associated to dilation on the odd fibers of~$\Pi TX$.

\begin{defn}
Let 
$$
\mathcal{O}(\Phi_0^{2|1}(X);\omega^{k/2})\subset \Gamma(\Phi_0^{2|1}(X);\omega^{k/2})
$$
denote the subspace of \emph{holomorphic sections} satisfying (\ref{eq:holo}). 
\end{defn}

\begin{proof}[Proof of Theorem~\ref{thethm}] 
The condition is equivalent to $F$ being holomorphic (i.e., $\partial_{\bar \ell}F=\partial_{\bar\ell'}F=0$), so restricting the map in Proposition~\ref{prop:Maass} to holomorphic Maass forms (i.e., weak modular forms) yields the isomorphism
$$
\bigoplus_{i+j=k} \Omega^j_{\rm cl}(X;\MF^i)\stackrel{\sim}{\to} \mathcal{O}(\Phi_0^{2|1}(X),\omega^{k/2}),
$$  
as claimed.
\ep

\subsection{Motivation from the $2|1$-dimensional sigma model}\label{sec:21motivate}

Physical intuition for $\Phi_0^{2|1}(X)$ comes from the $2|1$-dimensional sigma model with target~$X$. Let $S\times_\Lambda \R^2$ denote a family of Euclidean tori equipped with a fiberwise spin structure. An $S$-family of fields for the $2|1$-dimensional sigma model with target~$X$ consists of a map~$x\colon S\times_\Lambda \R^2 \to X$, and an odd spinor $\psi\in \Gamma(S\times_\Lambda \R^2,\Pi \overline K^{1/2} \otimes x^*TX)$ twisted by the pullback tangent bundle; here~$K$ denotes the canonical line on the family of tori. 
Choosing a Riemannian metric on~$X$, the classical action is
\beq\begin{array}{ccl}
\mathcal{S}(x,\psi)
&=&\displaystyle\frac{1}{2}\int_{\R^2/\Lambda} \left(\langle \partial_z x,\partial_{\bar z} x\rangle-i\langle \psi,\nabla_{z}\psi\rangle\right) d \bar z dz
\end{array}\label{eq:21classactcomp}
\eeq
where the integral is fiberwise over~$S$. In terms of the fields $(x,\psi)$, the Euler--Lagrange equations for this classical action (\cite{strings1} pages~654 and~664, or~\S\ref{sec:21Hess}) are
$$
\nabla_{\partial_{\bar z}} \partial_{z} x=\frac{1}{2} R(\psi,\psi)\partial_z x,\quad\quad \nabla_{\partial_z}\psi=0.
$$
The \emph{energy} of a family of classical solutions is the fiberwise integral of $\langle \partial_z x,\partial_{\bar z}x\rangle$ on~$S\times_\Lambda\R^2$. 

\begin{defn} A \emph{classical vacuum} over $S$ is an energy zero family of classical solutions. \end{defn}

\begin{lem} An $S$-family of fields $(x,\psi)$ is a classical vacuum if and only if~$x$ is a family of constant maps and~$\psi$ is a constant section. \end{lem}

\bp If the energy of an $S$-family of fields is zero, then the underlying maps of tori are constant, i.e., $x\colon S\times_\Lambda \R^2\to X$ factors through the projection $S\times_\Lambda \R^2\to S$ and is a family of constant maps. By virtue of being a classical solution, $\psi$ is anti-holomorphic, but since the tori are compact this requires $\psi$ be a constant section. \ep

We can repackage the data~$(x,\psi)$ in super space form as a map~$\phi$ from a family of super tori~$S\times_\Lambda\R^{2|1}$ to $X$. When $S\times_\Lambda\R^2$ is equipped with the odd (or nonbounding) fiberwise spin structures, the spinor bundle is trivializable and this family coincides with the super tori defining~$\Phi^{2|1}(X)$. In this case, an energy zero classical solution is a map $S\times_\Lambda \R^{2|1}\to X$ factoring through $S\times\R^{0|1}$: a trivialization of the spinor bundle identifies the covariantly constant section~$\psi$ with a map $S\times \R^{0|1}\to X$. 

Analogous to the $1|1$-dimensional case, we can promote the above discussion of fields as a sheaf of sets on the site of super manifolds to a sheaf of groupoids, i.e., a stack. We take as isomorphisms maps between super tori that preserve the action. In this case, these symmetries are the rigid conformal isometries between families of super tori. In summary, the stack $\Phi_0^{2|1}(X)$ has as objects over~$S$ classical vacua associated to odd spin structures on tori, and as morphisms over~$S$ the symmetries of the classical action. 

As in the previous case, we observe that although the classical action depends on a choice of metric on~$X$, the stack $\Phi_0^{2|1}(X)$ does not. 

%

\section{Warm-up 2: the perturbative $1|1$-sigma model and the $\hat{A}$-genus}\label{sec:warmup2}

Define the $\hat{A}$-class of a Riemannian manifold $X$ to be the polynomial in Pontryagin forms defined by the usual characteristic series~$\frac{x/2}{\sinh(x/2)}$, and denote this closed, even differential form by $\hat{A}(X)\in \Omega^{\rm ev}_{\rm cl}(X)$. By the proof of Proposition \ref{prop:Kthy}, we can view this element as a function on classical vacua, which we also denote by $\hat{A}(X)\in C^\infty (\Phi_0^{1|1}(X))$. The goal of the present section is to identify this function with the relative super determinant of a family of operators, denoted $\Delta^{1|1}_X$, parametrized by~$\Phi^{1|1}_0(X)$. 
The operators $\Delta^{1|1}_X$ depend on a choice of metric on~$X$, and are the kinetic operators encoding the Hessian of the classical action for supersymmetric mechanics restricted to the classical vacua; see~\S\ref{sec:11Hess}.  A similarly defined operator $\Delta^{1|1}_n$ serves to normalize our computations relative to a trivial bundle $\underline{\R}^n$ with $n={\rm dim}(X)$.

\begin{prop} For a Riemannian manifold~$X$, the relative $\zeta$-super determinant of the family of operators $\Delta^{1|1}_X$ gives a function 
$$
\frac{{\rm sdet}_\zeta(\Delta^{1|1}_X)}{{\rm sdet}_\zeta(\Delta^{1|1}_n)}=\hat{A}(X)\in C^\infty (\Phi^{1|1}_0(X))
$$
that agrees with the $\hat{A}$-form of~$X$ under the isomorphism of Proposition \ref{prop:Kthy}. 
For any 1-parameter family of Riemannian metrics on $X$, this construction determines a concordance $\hat{A}(X\times \R)\in C^\infty({\Phi}_0^{1|1}(X\times \R))$ between the representatives of the class $[\hat{A}(X)]$.

\label{prop:11zeta} 
\end{prop}

We will also formulate a variation that brings Bismut--Freed--Quillen determinant lines into the picture.
The pullback of $\Delta^{1|1}_X$ to any $S$-point of $\Phi_0^{1|1}(X)$ determines a trivial super determinant line bundle with metric, and for these trivializations to descend to the stack, we need a consistent choice of the square root of the norm squared of the determinant section. 


\begin{prop} 
The super determinant line bundles of $\Delta^{1|1}_X$ pulled back to objects of $\Phi_0^{1|1}(X)$ over~$S$ uniquely descend to a trivial line bundle on the stack~$\Phi_0^{1|1}(X)$ whose trivializing section~$\one$ satisfies $\|\one\|=1$ and the super determinant section in this trivialization is~$\hat{A}(X)\cdot \one.$ This metrized line bundle is natural with respect to isometries between Riemannian manifolds. 
\label{prop:Aroof}
\end{prop}

At the end of the section we will explain how the above constructions of $\hat{A}(X)$ define a local index map for $1|1$-Euclidean field theories that is closely related to the local index in the Atiyah--Singer theorem. 

The essence of the above pair of propositions has occurred in various guises elsewhere: most notably, we learned the basic pieces of the $\zeta$-determinant computation in our Proposition~\ref{prop:11zeta} from E.~Witten's article in \cite{strings1}, particularly pages~476-485. To some extent, our computations are a repackaging of Witten's core ideas in a different conceptual framework, e.g., designed to connect with the Segal--Stolz--Teichner program. Our approach is also very similar in spirit to that of R. Grady and O. Gwilliam in their construction of the $\hat{A}$-class~\cite{gradygwilliam}; there are some language barriers between the two approaches, but roughly their version of topological quantum mechanics arises from the large-volume limit of our classical field theory. 


\subsection{A family of operators on $\Phi_0^{1|1}(X)$}

The tangent space to $\Phi^{1|1}(X)$ is a vector bundle $T\Phi^{1|1}(X)$ whose fiber at an $S$-point $(r,\phi)$ has sections~$\Gamma(S\times_r \R^{1|1},\phi^*TX)$ as a $C^\infty(S)$-module; we pull this bundle back along families of isometries. Note $T\Phi^{1|1}(X)$ is an infinite-rank vector bundle for ${\rm dim}(X)>0$, as one should expect from a mapping space. The Riemannian metric and Levi-Civita connection on~$TX$ pull back to these spaces of sections, and together with the fiberwise volume form on $S\times_r \R$ we obtain a pairing on sections of~$T\Phi^{1|1}(X)$ at each $S$-point. Explicitly, for a section $\nu\in \Gamma(S\times_r \R^{1|1},\phi^*TX)$, and $i_0\colon S\times_r\R\hookrightarrow S\times_r\R^{1|1}$ the fiberwise inclusion of the reduced manifold, we define Taylor components
\beq
a:=i_0^*\nu,\quad\quad \eta=i_0^*((\phi^*\nabla)_{D}\nu)\label{eq:11Taylor}
\eeq
where $D=\partial_\theta-i\theta\partial_t$, and the inner product of sections $\nu$ and $\nu'$ is the sum of two terms
$$
\int_{S\times_r \R/S} \langle a,a'\rangle dt,\quad\quad  \int_{S\times_r \R/S} \langle \eta,\eta'\rangle dt
$$
where $\langle-,-\rangle$ denotes the pullback of the metric on~$TX$ to $S\times_r \R$. The pairing on the odd sections~$\eta$ comes from the pullback of the pairing on $\Gamma(\Pi TX)$,
\beq
(\Gamma(TX)\otimes \C^{0|1})\otimes (\Gamma(TX)\otimes \C^{0|1})&\stackrel{\sigma}{\cong}& (\Gamma(TX)\otimes \Gamma(TX))\otimes (\C^{0|1} \otimes \C^{0|1})\nonumber\\
&\cong& \Gamma(TX)\otimes \Gamma(TX)\stackrel{\langle-,-\rangle}{\longrightarrow} C^\infty(X)\nonumber
\eeq
where we have identified the module $\Gamma(\Pi TX)$ of sections with $\Gamma(TX)\otimes \C^{0|1}$ for $\C^{0|1}$ is the odd super line (as a super vector space), $\sigma$ is the braiding isomorphism, and we use the isomorphism of super vector spaces~$\C^{0|1}\otimes \C^{0|1}\cong \C$. 

We can restrict the vector bundle $T\Phi^{1|1}(X)$ to $\Phi_0^{1|1}(X)\subset \Phi^{1|1}(X)$, and the pairing picks out a subbundle that is the normal bundle to the inclusion of stacks. 

\begin{defn} Define $\mathcal{N}\Phi^{1|1}_0(X)\subset T\Phi^{1|1}(X)|_{\Phi_0^{1|1}(X)}$ as having $S$-points sections in the orthogonal complement of the constant sections, where a section $\nu$ is constant if~$\nabla_{\partial_t}\nu=0$. We use the notation $\Gamma_0(S\times_r \R^{1|1},\phi^*TX)\subset \Gamma(S\times_r \R^{1|1},\phi^*TX)$ to denote this orthogonal complement at an $S$-point $(r,\phi)$. 
\end{defn}

We will describe an exponential map on sections of $\mathcal{N}\Phi^{1|1}_0(X)$ in~\S\ref{sec:11exp}, so that these can be viewed as a tubular neighborhood of the substack~$\Phi_0^{1|1}(X)\subset \Phi^{1|1}(X)$. 

To simplify the notation, set $\nabla=\phi^*\nabla$. Define a function on sections of $\mathcal{N}\Phi^{1|1}_0(X)$ by
\beq
\Hess_\phi(\nu):=-i\int_{S\times_r\R^{1|1}/S} \langle \nu,\nabla_{\partial_t}\nabla_D\nu\rangle dtd\theta,\quad \nu\in \Gamma_0(S\times_r\R^{1|1},\phi^*TX)\label{eq:11Hess}
\eeq
where $D=\partial_\theta-i\theta\partial_t$ is the right-invariant vector field on the family $S\times_r \R^{1|1}$, and the integral is the Berezinian integral over the fibers of the projection $S\times_r \R^{1|1}\to S$. Since it is built out of right-invariant vector fields, the function $\Hess$ is automatically invariant under the left action by isometries. Therefore $\Hess$ defines a function on the stack as claimed. A component-form version of this function will facilitate computations. 

\begin{lem} 
Taylor expanding $\nu$ using (\ref{eq:11Taylor}) and performing the Berezin integral in (\ref{eq:11Hess}),
\beq
\Hess_\phi(\nu)&=&\Hess_\phi(a,\eta)=-\int_{S\times_r \R/S}\langle(\Delta^{1|1}_X)^{\ev} a,a\rangle+\langle(\Delta^{1|1}_X)^{\odd}\eta,\eta\rangle dt,\nonumber\\
&& (\Delta^{1|1}_X)^{\ev}:=-\nabla_{\partial_t}^2+\frac{i}{2} \cdot \mathcal{R}\nabla_{\partial_t}, \quad (\Delta^{1|1}_X)^{\odd}=i\cdot \nabla_{\partial_t}\nonumber
\eeq
where $\mathcal{R}=i_0^*\phi^*\mathcal{R}(D,D)$ is the $\End(\phi^*TX)$-valued function on $S\times_r\R$ determined by the curvature 2-form of the Levi--Civita connection. We use the notation $\Delta^{1|1}_X:=(\Delta_X^{1|1})^{\ev}\oplus (\Delta^{1|1}_X)^{\odd}$ to denote the operator direct sum. 
\end{lem}

\bp
We compute the Taylor components as in (\ref{eq:11Taylor}) of the section~$\nabla_{\partial_t}\nabla_D\nu$:
\beq
i^*_0(\nabla_{\partial_t}\nabla_D\nu)&=&\nabla_{\partial_t} i^*_0 (\nabla_D\nu)=\nabla_{\partial_t} \eta\nonumber\\
i^* _0\nabla_D(\nabla_{\partial_t}\nabla_D\nu)&=&i^*_0(\nabla_{\partial_t}\nabla_D\nabla_D \nu)=\frac{1}{2}i^*_0(\nabla_{\partial_t} (\mathcal{R}(D,D)-\nabla_{[D,D]})\nu)\nonumber\\
&=&(\mathcal{R}(D,D)\nabla_{\partial_t}+i\nabla_{\partial_t}^2)i^*_0\nu=(i\nabla_{\partial_t}^2+\mathcal{R}(D,D)\nabla_{\partial_t})a\nonumber
\eeq
where we used 
$$
\nabla_D^2=\frac{1}{2}(\nabla_D\nabla_D+\nabla_D\nabla_D)=\frac{1}{2}(\mathcal{R}(D,D)-\nabla_{[D,D]})=\frac{1}{2}\mathcal{R}(D,D)+i\nabla_{\partial_t}.
$$
So now we have
\beq
\Hess_\phi(\nu)&=&-i\int_{S\times_r\R^{1|1}/S}\langle a+\theta \eta,\nabla_{\partial_t} \eta+\theta(i\nabla_{\partial_t}^2+\frac{1}{2}\mathcal{R}\nabla_{\partial_t})a\rangle d\theta dt\nonumber\\
&=&-\int_{S\times_r\R} \langle a,(-\nabla_{\partial_t}^2+\frac{i}{2}\mathcal{R}\nabla_{\partial_t})a\rangle+\langle \eta,i\nabla_{\partial_t} \eta\rangle dt,\nonumber
\eeq
as claimed. 
\ep

We define a similar vector bundle over $\Phi_0^{1|1}(X)$ coming from a trivial bundle on $X$. 

\begin{defn} 
Define $\mathcal{Z}_n(\Phi_0^{1|1}(X))$ as having $S$-points sections $\Gamma(S\times_r \R^{1|1},\phi^*\underline{\R}^n)$ in the orthogonal complement of the constant sections, for $\underline{\R}^n$ the trivial bundle on $X$ with $n={\rm dim}(X)$. 
\end{defn}

Define a family of operators $\Delta_n^{1|1}$ coming from the function on sections
$$
\int_{S\times_r \R^{1|1}/S}\langle \nu,\partial_t D\nu\rangle,\quad\quad \nu\in \Gamma_0(S\times_r\R^{1|1},\phi^*\underline{\R}^n)
$$
and adapting the computations above to the trivial bundle with trivial connection, we have
$$
(\Delta^{1|1}_n)_{\rm bos}=-\partial_t^2,\quad\quad (\Delta^{1|1}_n)_{\rm fer}=i\partial_t. 
$$

\subsection{The $\hat{A}$-class as a relative $\zeta$-determinant}\label{sec:11kinetic}

\begin{proof}[Proof of Proposition \ref{prop:11zeta}]

We will identify the $\zeta$-super determinant of~$\Delta^{1|1}_X$ with the~$\hat{A}$-class in its form gotten from the characteristic series on the right hand side,
$$
\frac{x/2}{\sinh(x/2)}=\exp\left(\sum_{k=1}^\infty\frac{x^{2k}}{2k(2\pi i)^{2k}}2\zeta(2k)\right).
$$
To set up the computation, we pullback $\Delta_X^{1|1}$ and $\Delta_n^{1|1}$ along the map $\pi \colon \R_{>0}\times \Pi TX\to \Phi_0^{1|1}(X)$, producing families of operators acting on the bundles whose fiber at $r\in \R_{>0}$ is $C^\infty(\R/r\Z)\otimes \Gamma(\ev^*TX)$ for $\ev\colon \R^{0|1}\times \Pi TX\to X$ the evaluation and $C^\infty(\R/r\Z)\otimes \Gamma(\ev^*\underline{\R}^n)$, respectively. In this description of sections, $(\pi^*\phi^*\nabla)_{\partial_t}=d/dt\otimes \id_{TX}$ and so
$$
\pi^*(\Delta^{1|1}_X)^{\ev}=-\frac{d^2}{d t^2}\otimes \id_{TX}+i\frac{d}{dt} \otimes \mathcal{R},\quad\quad \pi^*(\Delta^{1|1}_X)^{\odd}=i\frac{d}{dt}\otimes \id_{TX},
$$
where now $\mathcal{R}$ is the $\End(p^*TX)$-valued function on $\Pi TX$ associated to the curvature 2-form. We use the basis for functions on $\R/r\Z$ given by~$F_n=e^{2\pi i nt/r}$, giving the $\zeta$-functions,
\beq
\zeta_X^{\ev}(s)&=&\sum_{n\ne 0} {\rm Tr}\left(\frac{4\pi^2n^2}{r^2}\otimes \id_{TX}+\frac{2\pi i n}{r}\otimes i\mathcal{R}\right)^s\nonumber\\
\zeta_X^{\odd}(s)&=&\sum_{n\ne 0} {\rm Tr}\left(-\frac{2\pi n}{r}\otimes \id_{TX}\right)^s\nonumber
\eeq
corresponding to the operators $\pi^*(\Delta^{1|1}_X)^{\ev}$ and $\pi^*(\Delta^{1|1}_X)^{\odd}$, respectively. Similarly, for~$\pi^*\Delta_n^{1|1}$ we get
$$
\zeta_n^{\ev}(s)=\sum_{n\ne 0} {\rm Tr}\left(\frac{4\pi^2n^2}{r^2}\otimes \id_{\underline{\R}^n}\right)^s,\quad\quad \zeta_n^{\odd}(s)=\sum_{n\ne 0} {\rm Tr}\left(-\frac{2\pi n}{r}\otimes \id_{\underline{\R}^n}\right)^s
$$
from the bosonic and fermionic parts. The contributions of $\zeta_n^{\odd}$ and $\zeta_{X}^{\rm odd}$ only depend on the dimension of the vector bundles $TX$ and $\underline{\R}^n$, and so their net contribution to relative super determinant is~1. Binomial expansion in odd variables gives
\beq
\zeta_X^{\ev}(s)&=&\sum_{n\ne 0} {\rm Tr}\left(\id-\mathcal{R}\otimes \frac{r}{2\pi n}\right)^s\left(\frac{4\pi^2 n^2}{r^2}\right)^s\nonumber\\
&=&\sum_{n\ne 0}\sum_{k=0}^{\rm finite} {\rm Tr}\left(\mathcal{R}^k\frac{s(s-1)\cdots (s-k+1)}{k!(2\pi n)^k} r^k\right)\left(\frac{4\pi^2 n^2}{r^2}\right)^s\nonumber
\eeq
where the sum over $k$ is finite because $\mathcal{R}$ is nilpotent. The $k=0$ contribution cancels identically with the contribution from $\zeta_n^{\ev}$. For $k>1$, we differentiate under the sum and obtain the contribution to $\zeta'(0)$
$$
\sum_{k=1}^{\rm finite} {\rm Tr}\left(\mathcal{R}^{k}\right)\frac{(-1)^{k-1}}{k(2\pi)^{k}} r^{k} 2\zeta(k)=-\sum_{k=1}^\infty\frac{{\rm Tr}(\mathcal{R}^{2k}) r^{2k}}{2k (2\pi)^{2k}}2\zeta(2k)
$$
where $\zeta(k)$ denotes the value of the Riemann $\zeta$-function at~$k$, and we have used that traces of odd powers of $\mathcal{R}$ vanish. So we have
$$
\frac{{\rm sdet}_\zeta(\Delta^{1|1}_X)}{{\rm sdet}_\zeta (\Delta^{1|1}_n)}=\exp\left(\sum_{k=1}^\infty\frac{{\rm Tr}(\mathcal{R}^{2k}) r^{2k}}{2k (2\pi i )^{2k}}\zeta(2k)\right)
$$
which by inspection defines an element of $C^\infty(\Phi_0^{1|1}(X))$; it remains to compare with the $\hat{A}$-form. 
In our cochain model, following considerations in subsection \ref{sec:chern} we have
$$
r^{2k}{\rm Tr}\left(\mathcal{R}\right)^{2k}=2(2k)! {\rm ph}_{k}(TX),
$$
where ${\rm ph}_{k}$ denotes the $4k^{\rm th}$ component of the Pontryagin character as a function on~$\Phi^{1|1}_0(X)$. Putting this together we get
\beq
\frac{{\rm sdet}_\zeta(\Delta^{1|1}_X)}{{\rm sdet}_\zeta (\Delta^{1|1}_n)}&=&\exp\left(\sum_{k=1}^\infty \frac{(2k)! {\rm ph}_{k}(TX)}{2k} \frac{2\zeta(2k)}{(2\pi i)^{2k}}\right)\nonumber
\eeq
which we identify as the $\hat{A}$-class on $X$ as a function on~$\Phi^{1|1}_0(X)$. 

Given a path of metrics on~$X$, we obtain a metric on $X\times \R$ via the direct sum of the standard metric on~$\R$ and the 1-parameter family of metric on~$X$. We may apply the construction above to~$X\times \R$ with this metric, yielding $\hat{A}(X\times \R)\in C^\infty(\Phi^{1|1}_0(X\times \R))$. Since the Pontryagin classes are stable, and $i_a^*T(X\times \R)\cong TX\oplus \underline\R$ for $i_a\colon X\hookrightarrow X\times \R$, the pullback of the relative $\zeta$-determinant at $a\in \R$ does indeed agree with the $\hat{A}$-form of~$X$ for the metric at~$a\in \R$. 
\ep


\begin{rmk} 
We will have some use below for the square root of the $\zeta$-determinant of $d^2/dt^2\otimes {\rm Id}_{TX}$: as a function on $\R_{>0}\times \Pi TX$ it is $r^{n}$, which follows from Example~2 of~\cite{zetaprod}.\label{rmk:11zeta}
\end{rmk} 

\subsection{The $\hat{A}$-class as a section of a Quillen determinant line}

\begin{proof}[Proof of Proposition~\ref{prop:Aroof}]

At each $S$-point of $\Phi^{1|1}_0(X)$, the pullback of the family of operators $\Delta^{1|1}_X$ and $\Delta^{1|1}_n$ defines a relative super determinant line bundle on $S$, which is the trivial line bundle whose trivializing section has norm squared the pullback of the function
\beq
\frac{{\rm sdet}_\zeta((\Delta^{1|1}_X)^*\Delta^{1|1}_X)}{{\rm sdet}_\zeta((\Delta^{1|1}_X)^*\Delta^{1|1}_n)}&=&\frac{\| {\rm det}_\zeta((d^2/dt^2\otimes {\rm Id}_{\underline{\R}^n})(d^2/dt^2\otimes {\rm Id}_{\underline{\R}^n})^*)\|^{1/2}}{\| {\rm det}_\zeta((d^2/dt^2\otimes {\rm Id}_{TX}+id/dt\mathcal{R})(d^2/dt^2\otimes {\rm Id}_{TX}+id/dt\mathcal{R})^*)\|^{1/2}}\nonumber\\
&=&\hat{A}(X)\cdot \hat{A}(X)\nonumber
\eeq
on $\Phi_0^{1|1}(X)$. A square root of the pullback of $\hat{A}(X)\cdot \hat{A}(X)$ defines a metric trivialization of this relative super determinant line over $S$. For these trivializations to descend to the stack, we need to make a universal choice of square root. The obvious one is~$\hat{A}(X)$ itself, which is uniquely characterized by being real and having the property at any $S$-point its restriction to the reduced manifold is the constant function~1. For this family of metric trivializations, the relative super determinant section is exactly~$\hat{A}(X)$.
\ep

\subsection{From the Hessian of the classical action to $\Delta^{1|1}_X$}\label{sec:11Hess}

For $\phi\colon S\times_r \R^{1|1}\to X$ a field, the superspace form of the classical action (\ref{eq:11classactcomp}) is
$$
\mathcal{S}(\phi)=\frac{i}{2}\int_{S\times_r\R^{1|1}} \langle \partial_t\phi,D\phi\rangle d\theta dt.
$$
Let $v$ and $w$ be sections of the tangent bundle at the point $(r,\phi)$, and let $\delta_v$, $\delta_w$ denote the derivations on functions associated to these vector fields. We compute the Hessian of $\mathcal{S}(\phi)$ on $\Phi_0^{1|1}(X)$, using that at a cricial point (i.e., where the first derivative vanishes) the Hessian can be computed as $\delta_w\delta_v\mathcal{S}$. 

The first derivative is 
\beq
\delta_v\mathcal{S}(\phi)&=&\frac{i}{2}\int_{S\times_r \R^{1|1}} \Big(\langle \nabla_v \partial_t \phi,D\phi\rangle+\langle \partial_t\phi,\nabla_v D\phi\rangle\Big)d\theta dt\nonumber\\
&=&\frac{i}{2}\int_{S\times_r \R^{1|1}} \Big(\langle \nabla_t \delta_v\phi,D\phi\rangle+\langle\partial_t\phi,\nabla_D \delta_v\phi\rangle \Big)d\theta dt\nonumber\\
&=&\frac{i}{2}\int_{S\times_r \R^{1|1}} \Big(-\langle \delta_v\phi,\nabla_t D\phi\rangle+\partial_t\langle \delta_v\phi,D\phi\rangle-\langle \nabla_D\partial_t\phi,\delta_v\phi\rangle+D\langle\partial_t\phi,\delta_v\phi\rangle\Big)d\theta dt\nonumber\\
&=&-i\int_{S\times_r \R^{1|1}}\langle \nabla_tD\phi,\delta_v\phi\rangle d\theta dt\nonumber
\eeq
where in the second line we use that the connection is torsion free, in the third line we integrate by parts, and in the fourth line we discarded terms that integrate to zero. This recovers the Euler--Lagrange equations in super space form, $\nabla_tD\phi=0.$

The 2nd order variation on the classical solutions (i.e., with $\nabla_tD\phi=0$) is
\beq
\delta_w\delta_vS(\phi)&=&-i\int_{S\times_r \R^{1|1}}\Big( \langle \nabla_w \delta_v\phi,\nabla_tD\phi\rangle+\langle \delta_v\phi,\nabla_w\nabla_tD \phi\rangle\Big)d\theta dt\nonumber\\
&=&-i\int_{S\times_r \R^{1|1}}\langle \delta_v\phi,\nabla_t\nabla_wD \phi\rangle d\theta dt\nonumber\\
&=&-i\int_{S\times_r \R^{1|1}}\langle \delta_v\phi,\nabla_t\nabla_D \delta_w \phi\rangle d\theta dt\nonumber
\eeq
where in the first line we enforce the Euler--Lagrange equation and $\nabla_t\nabla_w-\nabla_w\nabla_t=\mathcal{R}(\partial_t,w)=0$, and in the last line we use that the connection is torsion free. Hence, we have 
$$
\Hess(v,w)=\delta_v\delta_w S(\phi)=-i\int_{S\times_r \R^{1|1}}\langle v,\nabla_t\nabla_D w \rangle d\theta dt.
$$
Our construction of the $\hat{A}$-class uses this Hessian in a 1-loop quantization procedure as overviewed in~\S\ref{sec:introphys}.

\subsection{An exponential map for sections of $T\Phi^{1|1}(X)$ over $\Phi^{1|1}_0(X)$}\label{sec:11exp}

%

We now will describe an exponential map $T\Phi^{1|1}(X)|_{\Phi_0^{1|1}(X)}\to \Phi^{1|1}(X)$; this basically places E.~Witten's discussion on page 482 of \cite{strings1} in the context of our stacks of fields, and allows one to view sections of~$\mathcal{N}\Phi_0^{1|1}(X)$ as a tubular neighborhood of $\Phi_0^{1|1}(X)$ inside $\Phi^{1|1}(X)$.

Given a section $\nu\in \Gamma(S\times_r \R^{1|1},\phi^*TX)$ along a map $\phi$ that factors through the projection to $S\times \R^{0|1}$, we will use the Levi-Civita connection on $TX$ and the exponential map of the Riemannian manifold $X$ to define a map $S\times_r \R^{1|1}\to X$. We identify $\nu$ with a map $S\times \R^{1|1}\to TX$ satisfying a periodicity condition. Consequently, for any compact subset of~$S$ the image of~$\nu$ is a compact subset of~$TX$. We consider the composition
$$
[0,\delta]\times S\times \R^{1|1}\to TX\stackrel{\exp}{\to} X
$$
for $\exp$ the exponential map with respect to the Riemannian metric; there exists a fixed $\delta$ for which this exponential map is defined provided that $S$ is compact, or if one insists on noncompact families we may choose a smooth strictly positive function on $S$ we also denote by~$\delta$, and we consider the bundle of compact intervals $[0,\delta]\times S\subset \R\times S$ over~$S$. We denote the family of maps parametrized by $[0,\delta]\times S$ defined by the composition above by $\phi+\delta\nu$. 

To define the normal bundle to the inclusion, we need to identify the orthogonal complement to sections $\nu$ such that $\phi+\delta\nu\in \Phi^{1|1}_0(X)\subset \Phi^{1|1}(X)$, i.e., sections whose image under the exponential map remains in the substack $\Phi_0^{1|1}(X)$. Such sections are precisely constant ones, meaning those in the kernel of $(\phi^*\nabla)_{\partial_t}$. 

\subsection{A local index map for $1|1$-EFTs} \label{sec:11index}
When $X$ is oriented, the super manifold $\Pi TX$ has a canonical volume form coming from integration of differential forms on~$X$. Any nonvanishing function on $\Pi TX$ can be used to modify the volume form. Similarly, given a function on $\R_{>0}\times \Pi TX$ there is a canonical volume form along the fibers of the projection,~$\R_{>0}\times \Pi TX\to \R_{>0}$, that we can modify by any nonvanishing function on $\R_{>0}\times \Pi TX$. 

It follows from Remark~\ref{rmk:11zeta} that the $\zeta$-super determinant of $\Delta^{1|1}_X$ (\emph{not} the relative determinant) is~$r^{-n/2}\hat{A}(X)$. Together with the canonical volume form on the fibers of $\R_{>0}\times \Pi TX\to \R_{>0}$ we get a map on sections, 
$$
\Gamma (\Phi_0^{1|1}(X);\omega^{k/2})\stackrel{\cdot r^{-n/2}\hat{A}(X)}{\to} C^\infty(\R_{>0}\times \Pi TX)\stackrel{\int_X}{\to} \Gamma (\Phi_0^{1|1}(\pt);\omega^{k/2-n/2})
$$
where $\int_X$ denotes integration of differential forms over $X$ and ${\rm dim}(X)=n$. In this way, ${\rm sdet}_\zeta(\Delta^{1|1}_X)$ determines a pushforward along the map $\Phi_0^{1|1}(X)\to \Phi_0^{1|1}(\pt)$. 
The local index theorem (or the KO-version of Riemann--Roch) identifies this with the wrong-way map coming from the spin orientation of ${\rm KO}$ tensored with $\C$. Furthermore, the total volume of~$\Phi_0^{1|1}(X)$ with respect to this choice of volume form is the $\hat{A}$-genus of~$X$. 

We obtain an index map for $1|1$-EFTs over an oriented manifold $X$ by precomposing with the restriction map that evaluates a field theory on the super circles over $X$ defining the stack of classical vacua:
$$
1|1\EFT^k(X)\stackrel{\rm res}{\to} \Gamma (\Phi_0^{1|1}(X);\omega^{k/2})\stackrel{\int_X -\cdot r^{-n/2}\hat{A}(X)}{\longrightarrow} \Gamma (\Phi_0(\pt);\omega^{k/2-n/2})=\left\{ \begin{array}{cc} \C & k-n \ {\rm even} \\ 0 & k-n \ {\rm odd.}\end{array}\right.
$$
For field theories over $X$ defined in terms of vector bundles with connection, we showed in~\S\ref{sec:chern} that the first map gives ${\rm Tr}(\exp(irF))\in C^\infty (\Phi_0^{1|1}(X))$ for $F$ the curvature of the connection. Tracing through the various factors of $\pi $ and $r$, we see that its image under the above composition is the usual index map for the Atiyah--Singer theorem, namely we obtain the $\hat{A}$-genus of $X$ twisted by the Chern character of the bundle. 

\section{The perturbative $2|1$-sigma model and the Witten genus}



In this section we construct ${\rm Wit}^*(X)$ as a function on $\Phi_0^{2|1}(X)$ via the $\zeta$-determinant of a family of operators, denoted $\Delta_{X}^{2|1}$. This is a straightforward generalization of the $1|1$-dimensional case, and proves the first part of Theorems~\ref{thm:Witzeta} and~\ref{thm2}. The main new feature of the $2|1$-dimensional situation is the rational string obstruction: ${\rm Wit}^*(X)$ only determines a cocycle if the first Pontryagin \emph{form} happens to vanish. When this form doesn't vanish, a choice rational string structure specifies a concordance from~${\rm Wit}^*(X)$ to a cocycle. 
This is the second part of Theorems~\ref{thm:Witzeta} and~\ref{thm2}.


\subsection{A family of operators on $\Phi_0^{2|1}(X)$}

The tangent space to $\Phi^{2|1}(X)$ is a vector bundle $T\Phi^{2|1}(X)$ whose fiber at an $S$-point $(r,\phi)$ is~$\Gamma(S\times_\Lambda \R^{2|1},\phi^*TX)$ as (an infinite-rank) $C^\infty(S)$-module. We pull this module back along families of isometries, yielding a vector bundle over the stack. The Riemannian metric and Levi-Civita connection on~$TX$ pull back to these spaces of sections, and the fiberwise volume form on $S\times_\Lambda \R^2$ gives a pairing on sections of~$T\Phi^{2|1}(X)$ at each $S$-point. This is completely analogous as in the $1|1$-dimensional case, where we define Taylor components
\beq
a:=i_0^*\nu,\quad\quad \eta=i_0^*(\phi^*\nabla)_D\nu\label{eq:21Taylor}
\eeq
for $D=\partial_\theta+\theta\partial_{\bar z}$ and $i_0\colon S\times_\Lambda\R^2\hookrightarrow S\times_\Lambda \R^{2|1}$ the inclusion along the fiberwise reduced manifold. Then we have pairings
$$
\int_{S\times_\Lambda \R^2/S} \langle a,a'\rangle \frac{1}{2i}d\bar z dz,\quad\quad \int_{S\times_\Lambda \R^2/S} \langle \eta,\eta'\rangle \frac{1}{2i}d\bar z dz.
$$

We can restrict $T\Phi^{2|1}(X)$ to $\Phi_0^{2|1}(X)\subset \Phi^{2|1}(X)$, and the pairing picks out a subbundle that is the normal bundle to the inclusion of stacks. 

\begin{defn} Define $\mathcal{N}\Phi^{2|1}_0(X)\subset T\Phi^{2|1}(X)|_{\Phi_0^{2|1}(X)}$ as having $S$-points sections $\Gamma(S\times_\Lambda \R^{2|1},\phi^*TX)$ in the orthogonal complement of the constant sections. We use the notation $\Gamma_0(S\times_\Lambda \R^{2|1},\phi^*TX)\subset \Gamma(S\times_\Lambda\R^{2|1},\phi^*TX)$ to denote this orthogonal complement at an $S$-point $(\Lambda,\phi)$. 
\end{defn}

The exponential map in~\S\ref{sec:11exp} for the $1|1$-dimensional case continues to make sense here, and allows one to view sections of this normal bundle as defining a tubular neighborhood of $\Phi_0^{2|1}(X)$ inside $\Phi^{2|1}(X)$. To normalize various computations, we define a similar vector bundle over $\Phi_0^{2|1}(X)$ coming from a trivial bundle on $X$. 

\begin{defn} 
Define $\mathcal{Z}_n(\Phi_0^{2|1}(X))$ as having $S$-points sections $\Gamma(S\times_\Lambda \R^{2|1},\phi^*\underline{\R}^n)$ in the orthogonal complement of the constant sections, for $\underline{\R}^n$ the trivial bundle with trivial connection on $X$ and $n={\rm dim}(X)$. 
\end{defn}

To simplify the notation, set $\nabla=\phi^*\nabla$. Define a function on sections of $\mathcal{N}\Phi^{2|1}_0(X)$ by
\beq
\Hess_\phi(\nu):=\int_{S\times_\Lambda\R^{2|1}/S} \langle \nu,\nabla_{\partial_z}\nabla_D\nu\rangle d\bar z d zd\theta,\quad \nu\in \Gamma_0(S\times_\Lambda\R^{2|1},\phi^*TX)\label{eq:21Hess}
\eeq
where $D=\partial_\theta+\theta\partial_{\bar z}$ is the right-invariant vector field on the family $S\times \R^{2|1}$, and the integral is the Berezinian integral over the fibers of the projection $S\times_\Lambda \R^{2|1}\to S$. Since it is built out of right-invariant vector fields on~$\R^{2|1}$, the function~$\Hess$ is automatically invariant under the left action of isometries. Therefore it defines a function on the stack as claimed. A component-form version of this function will facilitate computations. 

\begin{lem} 
Taylor expanding $\nu$ using (\ref{eq:21Taylor}) and performing the Berezin integral in (\ref{eq:11Hess}),
\beq
\Hess_\phi(\nu)&=&\Hess_\phi(a,\eta)=\int_{S\times_\Lambda \R^2/S}\langle(\Delta^{2|1}_X)^{\ev} a,a\rangle+\langle(\Delta^{2|1}_X)^{\odd}\eta,\eta\rangle d\bar z dz,\nonumber\\
&& (\Delta^{2|1}_X)^{\ev}:=-\nabla_{\partial_z}\nabla_{\partial_{\bar z}}+\frac{1}{2} \mathcal{R}\nabla_{\partial_z}, \quad (\Delta^{2|1}_X)^{\odd}= \nabla_{\partial_z}\nonumber
\eeq
where $\mathcal{R}:=i_0^*\phi^*\mathcal{R}(D,D)$ is the $\End(\phi^*TX)$-valued function on $S\times_\Lambda\R^{2}$ determined by the curvature 2-form of the Levi--Civita connection.
\end{lem}

The proof is identical to the $1|1$-dimensional case. Let $\Delta_X^{2|1}:=(\Delta_X^{2|1})^{\ev}\oplus (\Delta^{2|1}_X)^{\odd}$ denote the operator direct sum. This defines the relevant family of operators over $\Phi_0^{2|1}(X)$. 


\subsection{The non-holomorphic Witten class as a $\zeta$-determinant}\label{prob1}

\begin{proof}[Proof of Theorem \ref{thm:Witzeta}] 

We pullback $\Delta_X^{2|1}$ along the map $\pi \colon L\times \Pi TX\to \Phi_0^{2|1}(X)$, producing a family of operators acting on the bundle whose sections at $\Lambda \in L$ is $C^\infty(\R^{2}/\Lambda)\otimes \Gamma(\ev^*TX)$ for $\ev\colon \R^{0|1}\times \Pi TX\to X$ the evaluation map. In this description of sections, $(\pi^*\phi^*\nabla)_{\partial_z}=\partial_z\otimes \id_{TX}$ and $(\pi^*\phi^*\nabla)_{\partial_{\bar z}}=\partial_{\bar z}\otimes \id_{TX}$, so
$$
\pi^*(\Delta^{2|1}_X)^{\ev}=-\partial_z\partial_{\bar z}\otimes \id_{TX}+\partial_z \otimes \mathcal{R},\quad\quad \pi^*(\Delta^{2|1}_X)^{\odd}= \partial_z\otimes \id_{TX},
$$
where now $\mathcal{R}$ is the $\End(p^*TX)$-valued function on $\Pi TX$ associated to the curvature 2-form. We use the basis for functions on $\R^2/\ell\Z\oplus \ell'\Z$
\beq
F_{n,m}(z,\bar{z}):=\exp\left(\frac{\pi}{\vol}(-z(n \bar\ell +m\bar{\ell}')+\bar{z}(n\ell+m\ell'))\right),\quad (m,n)\in \Z\times \Z,\label{eq:funontori}
\eeq
where $\vol=(\ell\bar\ell'-\bar\ell\ell')/2i$.  We form the $\zeta$-functions,
\beq
\zeta_X^{\ev}(s)&=&\sum_{(m,n)\in \Z^2_*} {\rm Tr}\left(\frac{\pi^2}{\vol^2}|m\ell+n\ell'|^2\otimes {\rm Id}_{TX}+\frac{\pi}{\vol}(m\bar\ell+n\bar\ell')\otimes \mathcal{R}\right)^s\nonumber\\
\zeta_X^{\odd}(s)&=&\sum_{(m,n)\in \Z^2_*} {\rm Tr}\left(\frac{\pi}{\vol} (m\bar\ell+n\bar\ell')\otimes {\rm Id}_{TX}\right)^s\nonumber
\eeq
corresponding to the operators $\pi^*(\Delta_X^{2|1})^{\ev}$ and $\pi^*(\Delta_X^{2|1})^{\odd}$, respectively. Similarly for $\Delta^{2|1}_n$, we have
\beq
\zeta_n^{\ev}(s)&=&\sum_{(m,n)\in \Z^2_*} {\rm Tr}\left(\frac{\pi^2}{\vol^2}|m\ell+n\ell'|^2\otimes {\rm Id}_{\underline{\R}^n}\right)^s\nonumber\\
\zeta_n^{\odd}(s)&=&\sum_{(m,n)\in \Z^2_*} {\rm Tr}\left(\frac{\pi}{\vol} (m\bar\ell+n\bar\ell')\otimes {\rm Id}_{\underline{\R}^n}\right)^s\nonumber
\eeq
The Pfaffian contributions from operators acting on odd sections can be computed following Example~5 in~\cite{zetaprod}. However, since the $\zeta$-functions associated to these pfaffians are equal, their overall contribution cancels. For the operators on even sections, we take the binomial expansion,
\beq
\zeta_X^{\ev}(s)&=&\sum_{(m,n)\in \Z^2_*} {\rm Tr}\left(\left({\rm Id}_{TX}+\frac{\vol}{2\pi} (m\ell+n\ell')^{-1}\otimes\mathcal{R}\right)^s\left(\frac{\pi^2}{\vol^2}|m\ell+n\ell'|^2\right)^{s}\right) \nonumber\\
&=&\sum_{(m,n)\in \Z^2_*} \sum_{k=0}^{\rm finite} {\rm Tr}\left[\left(\vol^{k}\frac{s(s-1)\cdots (s-k+1)}{k!(2\pi)^{k}(m\ell+n\ell')^{k}}\otimes \mathcal{R}^k\right)\left(\frac{\pi |m\ell+n\ell'|}{\vol} \right)^{2s}\right]\nonumber
\eeq
where the sum is finite because~$\mathcal{R}$ is nilpotent. Focusing on the part of the sum starting from $k=3$, we differentiate under the sum to obtain the following contribution to $\zeta'(0)$:
$$
\sum_{(n,m)\in \Z^2_*} \sum_{k=3}^{\rm finite} {\rm Tr}\left(\frac{(-1)^{k-1}}{k} \frac{\vol^k}{(2\pi)^k} (m\ell+n\ell')^{-k}\otimes \mathcal{R}^{k}\right)=-\sum_{k=2}^\infty\frac{\vol^{2k} E_{2k}}{2k(2\pi)^{2k} }{\rm Tr}(\mathcal{R}^{2k})
$$
where we have used that odd powers of $\mathcal{R}$ have trace zero. For $k=0$, we obtain the same $\zeta$-function in the calculation of ${\rm det}_\zeta(\partial_z\partial_{\bar z}\otimes {\rm Id}_{\underline{\R}^n})$, so this cancels in the relative determinant. The derivative at $s=0$ of the $k=2$ term is essentially the definition of the non-holomorphic 2nd Eisenstein series, 
$$
\lim_{s\to 0^-} \frac{d}{ds}{\rm Tr}\left(\mathcal{R}^2\sum_{(m,n)\in \Z^2_*} \frac{s(s-1)}{2}\left(\frac{\vol^2 |m\ell+n\ell'|^{4s}(m\ell+n\ell')^{-2}}{(2\pi)^2}\right)\right)=-\frac{\vol^2E_2^*}{2(2\pi)^2}{\rm Tr}(\mathcal{R}^2).
$$
So altogether we have
$$
\frac{{\rm sdet}_\zeta\left(\Delta^{2|1}_X\right)}{{\rm sdet}_\zeta \left(\Delta^{2|1}_n\right)}=\exp\left(\frac{\vol^2{\rm Tr}(\mathcal{R}^2)}{4(2\pi i)^2}E_2^*+\sum_{k\ge 2}^\infty \frac{\vol^{2k} {\rm Tr}(\mathcal{R}^{2k})}{4k(2\pi i)^{2k}} E_{2k}\right),
$$
and next we identify this function on $\Phi^{2|1}_0(X)$ with the non-holomorphic Witten class. 

The inclusion
$$
\Omega^\bullet_{\rm cl}(X;\C)\cong \Omega^\bullet_{\rm cl}(X;\MF^0)\hookrightarrow \bigoplus_{i+j=\bullet} \Omega^i(X;\MF^j)\cong \Gamma (\Phi^{2|1}(X);\omega^{\bullet/2})
$$
specifies the Pontryagin character of $X$ in our cochain model,
\beq
2(2k)! {\rm ph}_{k}(TX)=\vol^{2k}{\rm Tr}(\mathcal{R}^{2k})\in \Gamma (\Phi_0^{2|1}(X);\omega^{\bullet/2}).\label{eq:TMFChern}
\eeq
So we have
$$
\frac{{\rm sdet}_\zeta\left(\Delta^{2|1}_X\right)}{{\rm sdet}_\zeta \left(\Delta^{2|1}_n\right)}=\exp\left(\frac{{\rm ph}_1(TX)}{(2\pi i)^2}E_2^*+\sum_{k\ge 2}^\infty \frac{(2k)!{\rm ph}_k(TX)}{2k(2\pi i)^{2k}} E_{2k}\right),
$$
which is the non-holomorphic Witten class as a smooth function on $\Phi_0^{2|1}(X)$, which by Proposition~\ref{prop:Maass} can be identified with an element of the ring of closed differential forms with values in weak Maass forms. 

A 1-parameter family of metrics on~$X$ determines a metric on $X\times \R$ via the direct sum of the standard metric on~$\R$ and the 1-parameter family of metric on~$X$. Applying the construction above to~$X\times \R$ yields ${\rm Wit}^*(X\times \R)\in C^\infty(\Phi^{2|1}_0(X\times \R))$. Since the Pontryagin classes are stable, and $i_a^*T(X\times \R)\cong TX\oplus \underline\R$ for $i_a\colon X\hookrightarrow X\times \R$, the pullback of the relative $\zeta$-determinant at $a\in \R$ does indeed agree with the Witten form of~$X$ for the metric at~$a\in \R$.

Finally, a rational string structure $H$ with $dH=-p_1(TX)$ specifies the function 
$$
\frac{{\rm sdet}_\zeta\left(\Delta^{2|1}_X\right)}{{\rm sdet}_\zeta \left(\Delta^{2|1}_n\right)}\cdot \exp\left( \frac{d(tH) }{(2\pi i)^2}E_2^* \right)\in C^\infty(\Phi_0^{2|1}(X\times \R)),
$$
that at $t=0$ is the original relative super determinant, and at $t=1$ we obtain the modular and holomorphic Witten class of~$X$. The concordance determined by~$H$ from $p_1$ to zero exists if and only if $X$ is rationally string, since (by Stokes' Theorem) a pair of closed differential forms are concordant if and only if they are cohomologous. 
\ep

\subsection{The super determinant line bundle over $\Phi_0^{2|1}(X)$}\label{sec:sdetline}

A different perspective on the string obstruction and the Witten class uses a mild generalization of the Bismut--Freed--Quillen formalism of determinant line bundles, as discussed in~\S\ref{sec:sdet}. For each $S$-point of $\Phi_0^{2|1}(X)$, we can pullback the family of operators $\Delta^{2|1}_X$ to $S$, yielding a family of invertible operators over~$S$. Define the \emph{relative super determinant line bundle} over~$S$ as the trivial line bundle whose trivializing section has norm squared the pullback of~$\det_\zeta((\Delta^{2|1}_X)^*\Delta^{2|1}_X)$. We call this trivializing section the \emph{relative super determinant section.}

To assemble these line bundles and sections at each $S$-point into a line bundle on the stack $\Phi_0^{2|1}(X)$, we require isomorphisms of metrized line bundles for each isomorphism between $S$-points. This basically amounts to fixing the phase of the determinant. Since all the relevant bundles over~$S$ are (topologically) trivial, it suffices to fix compatible metric trivializations at each $S$-point. There turns out to be an essentially unique way to do this, characterized by Lemma~\ref{lem:uniquefun}. 


\begin{proof}[Proof of Lemma~\ref{lem:uniquefun}] Since the adjoint of an operator has conjugate eigenvalues, a computation analogous to the one in the previous subsection shows that
$$
\frac{{\rm sdet}_\zeta((\Delta^{2|1}_X)^*\Delta^{2|1}_X)}{{\rm sdet}_\zeta((\Delta^{2|1}_n)^*\Delta^{2|1}_n)}={\rm Wit}^*(X)\cdot \overline{{\rm Wit}^*(X)}=\|{\rm Wit}^*(X)\|^2.
$$
Using the cover~$L\times \Pi TX\to \Phi_0^{2|1}(X)$, the universal restriction to the reduced manifold of~$S$ comes from the inclusion $L\times X\subset L\times \Pi TX$. We observe that indeed`${\rm Wit}^*(X)$ pulls back to~1 on this restriction. Hence, ${\rm Wit}^*(X)$ satisfies the requisite properties. 

Now let $\sigma$ be a second function satisfying the conditions in the lemma. Then ${\rm Wit}^*(X)=f\cdot \sigma$ for $f$ a $U(1)$-valued function on~$\Phi_0^{2|1}(X)$. Pulling $f$ back along the cover $p\colon L\times \Pi TX\to \Phi_0^{2|1}(X)$ we find
$$
p^*f=\exp(i\sum_k f_k\otimes \omega_k)
$$
where $f_k$ is a real-valued function on $L$ and $\omega_k$ is a real differential form on~$X$. For $f_k\otimes \omega_k$ to descend to a function on the stack, $f_k$ must be dilation-equivariant for the $\C^\times$-action, commensurate with the degree of $\omega_k$. But for $f_k$ to also be real-valued, the dilation action on $f_k$ must be trivial and hence the degree of $\omega_k$ must also be zero. But then condition~(2) in the lemma shows that $p^*f\equiv 1$, so $f\equiv 1$, proving the lemma. 
\ep

Rescaling the relative super determinant line at each $S$-point by the pullback of~$1/{\rm Wit}^*(X)$ gives a metric trivialization, and these trivializations are compatible since ${\rm Wit}^*(X)$ pulls back from the stack. This justifies the following definition.

\begin{defn} The \emph{relative super determinant line} of~$\Delta_X^{2|1}$ and its \emph{relative determinant section} is the trivial metrized line with section~${\rm Wit}^*(X)$. 
\end{defn}

The failure of ${\rm Wit}^*(X)$ to be a cocycle comes from it not satisfying equations (\ref{eq:holo}), which in turn stems from $E_2^*$ not being holomorphic. To measure this failure, we write
\beq
{\rm Wit}^*(X)&=&\exp\left(\frac{\vol^2{\rm Tr}(\mathcal{R}^2)}{2(2\pi i)^2}E_2^*+\sum_{k\ge 2}^\infty \frac{\vol^{2k} {\rm Tr}(\mathcal{R}^{2k})}{2k(2\pi i)^{2k}} E_{2k}\right),\nonumber\\
&=&\exp\left(\frac{\vol^2{\rm Tr}(\mathcal{R}^2)}{2(2\pi i)^2}(E_2-\frac{\pi}{\vol})+\sum_{k\ge 2}^\infty \frac{\vol^{2k} {\rm Tr}(\mathcal{R}^{2k})}{2k(2\pi i)^{2k}} E_{2k}\right)\nonumber\\
&=&{\rm Wit}(X)\cdot \exp\left(\frac{\vol \cdot{\rm Tr}(\mathcal{R}^2)}{8\pi }\right)\nonumber\\
&=&{\rm Wit}(X)\cdot \exp\left(\frac{p_1(TX)}{2\pi \cdot \vol}\right)\nonumber
\eeq
where ${\rm Wit}(X)$ is a function on $L\times \Pi TX$ satisfying~(\ref{eq:holo}). 
Hence, the canonical modification 
$$
\exp(-p_1(TX)/(2\pi \vol))\cdot {\rm Wit}^*(X)={\rm Wit}(X),
$$ 
mediates between the non-holomorphic and the non-modular Witten class. However, it is easy to see that $\exp(-p_1(TX)/(2\pi \vol))$ is not a function on $\Phi_0^{2|1}(X)$: it is not invariant under the $\C^\times$-action. Instead it is a section of the line bundle $\mathcal{S}tr$ that is concordant to the trivial line bundle if and only if~$X$ has a rational string structure. 

\begin{proof}[Proof of Theorem~\ref{thm2}] A concordance between $\mathcal Str$ and the trivial line is equivalent to a concordance between $\exp(-p_1(TX)/8\pi i \vol)$ and the constant function~$1$. In turn, this requires a concordance from $p_1(TX)$ to zero, which exists if and only if $X$ has a rational string structure. With a choice of rational string structure fixed, consider the function
$$
\tilde{\sigma}:=\exp\left(\frac{\vol^2d(\lambda H)}{(2\pi i)^2}E_2+\sum_{k\ge 2}^\infty \frac{\vol^{2k} {\rm Tr}(\mathcal{R}^{2k})}{2k(2\pi i)^{2k}} E_{2k}\right)
$$
on $L\times \Pi T(X\times \R)$. The transformation properties of this function define the line bundle $\widetilde{\mathcal Str}$ on $\Phi_0^{2|1}(X\times \R)$, that indeed restricts to $\mathcal Str$ at $X\times\{1\}$ and the trivial bundle at $X\times\{0\}$. Furthermore, the restriction of $\tilde\sigma$ agrees with ${\rm Wit}(X)$ at~$1$, and is ${\rm Wit}_H(X)\in \mathcal{O} (\Phi_0^{2|1}(X))$ when~$\lambda=0$. 
\ep

\subsection{From the Hessian of the classical action to $\Delta^{2|1}_X$}\label{sec:21Hess}

The computation of the Hessian is very similar to the $1|1$-dimensional case, but we include it for completeness. The superspace form of the classical action~(\ref{eq:21classactcomp}) is
$$
\mathcal{S}(\phi)=i\int_{S\times_\Lambda \R^{2|1}} \langle \partial_z\phi,D\phi\rangle d\theta d\vol=i\int_{S\times_\Lambda \R^{2|1}} \langle \partial_z\phi,D\phi\rangle d\theta \Big(\frac{i}{2}d\bar zd z\Big)
$$
Let $v$ and $w$ be sections of the tangent bundle at the point $(\Lambda,\phi)$, and let $\delta_v$, $\delta_w$ denote the action of these vector field as variations of the field~$\phi$. We compute the Hessian of $\mathcal{S}(\phi)$ for $(\Lambda,\phi)$ and $S$-point of $\Phi_0^{2|1}(X)$. As usual, at a critical point of $S$ (i.e., where the first order derivative vanishes) the Hessian can be computed as $\delta_w\delta_vS$. 

The first derivative is 
\beq
\delta_v\mathcal{S}(\phi)&=&-\frac{1}{2}\int_{S\times_\Lambda \R^{2|1}} \Big(\langle \nabla_v \partial_z \phi,D\phi\rangle+\langle \partial_z\phi,\nabla_v D\phi\rangle\Big)d\theta d\bar z dz\nonumber\\
&=&-\frac{1}{2}\int_{S\times_\Lambda \R^{2|1}} \Big(\langle \nabla_{\partial_z} \delta_v\phi,D\phi\rangle+\langle\partial_z\phi,\nabla_D \delta_v\phi\rangle \Big)d\theta d\bar z dz\nonumber\\
&=&-\frac{1}{2}\int_{S\times_\Lambda \R^{2|1}} \Big(-\langle \delta_v\phi,\nabla_{\partial_z} D\phi\rangle+\partial_z\langle \delta_v\phi,D\phi\rangle-\langle \nabla_D\partial_z\phi,\delta_v\phi\rangle+D\langle\partial_z\phi,\delta_v\phi\rangle\Big)d\theta d\bar z dz\nonumber\\
&=&\int_{S\times_\Lambda \R^{2|1}}\langle \nabla_{\partial_z}D\phi,\delta_v\phi\rangle d\theta d\bar z dz\nonumber
\eeq
where in the second line we use that the connection is torsion free, in the third line we integrate by parts, and in the fourth line we discarded terms that integrate to zero: 
\beq
\int_{S\times_\Lambda \R^{2|1}}\Big(\partial_z\langle \delta_v\phi,D\phi\rangle+D\langle\partial_z\phi,\delta_v\phi\rangle\Big)&=&\int_{S\times_\Lambda \R^{2|1}}\Big(\partial_z\langle \delta_vx+\theta\nabla_v \psi,\theta\partial_{\bar z} x+\psi\rangle \nonumber \\
&&{}+(\partial_\theta+\theta\partial_{\bar z})\langle \partial_z x+\theta\nabla_{\partial_z}\psi,\delta_vx+\theta\nabla_v\psi\rangle\Big)d\theta d\bar z dz\nonumber\\
&=&\int_{S\times_\Lambda\R^2}\Big(\partial_z\langle \delta_vx,\partial_{\bar z} x\rangle+\partial_z\langle \nabla_v\psi,\psi\rangle\nonumber\\
&&\partial_{\bar z}\langle \partial_zx,\delta_vx\rangle-\langle \nabla_{z}\psi,\nabla_v\psi\rangle+\langle \nabla_{z}\psi,\nabla_v\psi\rangle\Big)d\bar z dz\nonumber \\
&=&\int_{S\times_\Lambda \R^2} d\Big(\langle \partial_z x,\delta_vx\rangle dz-(\langle \delta_vx,\partial_{\bar z} x\rangle+\langle \nabla_v\psi,\psi\rangle) d\bar z \Big)\nonumber
\eeq
where $d=dz\partial_z+d\bar z\partial_{\bar z}$. This recovers the Euler--Lagrange equations in super space form,~$\nabla_{\partial_z}D\phi=0.$

The 2nd order variation on classical vacua is
\beq
\delta_w\delta_v\mathcal{S}(\phi)&=&\int_{S\times_\Lambda \R^{2|1}}\Big( \langle \nabla_w \delta_v\phi,\nabla_{\partial_z}D\phi\rangle+\langle \delta_v\phi,\nabla_w\nabla_{\partial_z}D \phi\rangle\Big)d\theta d\bar z dz\nonumber\\
&=&\int_{S\times_\Lambda \R^{2|1}}\langle \delta_v\phi,\nabla_{\partial_z}\nabla_wD \phi\rangle d\theta d\bar z dz\nonumber\\
&=&\int_{S\times_\Lambda \R^{2|1}}\langle \delta_v\phi,\nabla_{\partial_z}\nabla_D \delta_w \phi\rangle d\theta d\bar z dz\nonumber
\eeq
where in the first line we enforce the Euler--Lagrange equation and $\nabla_{\partial_z}\nabla_w-\nabla_w\nabla_{\partial_z}=\phi^*\mathcal{R}(\partial_z,w)=0$, and in the last line we use that the connection is torsion free. Hence, on the tangent space we have 
$$
\Hess(v,w)=\delta_v\delta_w \mathcal{S}(\phi)=\int_{S\times_\Lambda \R^{2|1}}\langle v,\nabla_t\nabla_D w \rangle d\theta d\bar z dz.
$$

\subsection{An index map for $2|1$-EFTs} \label{subsec:21index}
Just as in the $1|1$-dimensional case, there is a canonical volume form on the fibers of $\Phi^{2|1}_0(X)\to \Phi^{2|1}_0(\pt)$, and the Witten class can be used to modify the associated integration map. With ${\rm dim}(X)=n$, the rescaling by $\vol^{-n/2} {\rm Wit}_H(X)$ gives a pushforward on sections, 
\beq
\Gamma (\Phi_0^{2|1}(X);\omega^{k/2})\stackrel{\cdot \vol^{-n/2}{\rm Wit}_H(X)}{\to} C^\infty(L \times \Pi TX)\stackrel{\int_X}{\to} \Gamma (\Phi_0^{2|1}(\pt);\omega^{k-n/2})\label{eq:21pushforwardfinal}
\eeq
that preserves the subspace of cocycles (i.e., the holomorphic functions). By the local index theorem, (\ref{eq:21pushforwardfinal}) can be identified with the wrong-way map coming from the string orientation of ${\rm TMF}$ tensored with $\C$. Furthermore, the total relative volume of $\Phi_0^{2|1}(X)$ with respect to this choice of volume form is the Witten genus of~$X$. Notice that this type of volume form fails to exist when~$X$ is not rationally string. 

We obtain an index map for $2|1$-EFTs over a rational string manifold~$X$ by precomposing with the map (\ref{eq:CH}) that evaluates a field theory on the super tori over~$X$, 
$$
2|1\EFT^\bullet(X)\stackrel{Z}{\to} \Gamma (\Phi^{2|1}_0(X);\omega^{\bullet/2})\to  \Gamma (\Phi^{2|1}_0(\pt);\omega^{\bullet-n/2})\cong {\rm MF}^{\bullet-n}.
$$
Our computation in dimension $1|1$ equated the analogous map with the index of twisted Dirac operators. As such, we view the above as a candidate geometric construction of a local analytical index for~TMF. 

\bibliographystyle{amsalpha}
\bibliography{references}

\end{document}